\documentclass[english]{article}
\usepackage[T1]{fontenc}
\usepackage[latin9]{inputenc}
\usepackage{morefloats}
\usepackage[letterpaper]{geometry}
\usepackage{float}
\usepackage{amsmath}

\usepackage{graphicx}
\usepackage{color}
\usepackage{epstopdf}
\usepackage{adjustbox}
\usepackage{diagbox}
\usepackage{multirow}
\usepackage{array}
\usepackage{mathrsfs}
\usepackage{amsmath, amssymb, amsthm}
\usepackage{subfigure}
\usepackage{tikz}

\graphicspath{{./Figures/}}

\newtheorem{theorem}{Theorem}

\newtheorem{lemma}{Lemma}

\usepackage{amsfonts} 
\usepackage{color}
\usepackage{adjustbox}
\usepackage{diagbox}
\definecolor{black}{rgb}{0,0,0}

\definecolor{red}{rgb}{1,0,0}

\definecolor{blue}{rgb}{0,0,1}

\usepackage{multirow}

\makeatother

\usepackage{babel}
\usepackage{authblk}
\title{}

\title{\textbf{A comparison of mixed multiscale finite element methods for multiphase transport in highly heterogeneous media}}
\author{Yiran Wang\thanks{Department of Mathematics, The Chinese University of Hong Kong, Shatin, Hong Kong.}, \;
	Eric Chung\thanks{Department of Mathematics, The Chinese University of Hong Kong, Shatin, Hong Kong.} \;
	and \; Shubin Fu\thanks{Department of Mathematics, University of Wisconsin-Madison,WI, USA. Corresponding author (shubinfu89@gmail.com)}
}

\begin{document}
	\maketitle
	\begin{abstract}
		
	\end{abstract}
	In this paper, we systemically review and compare two mixed multiscale finite element methods (MMsFEM) for multiphase transport in highly heterogeneous media. In particular, we will consider
	the mixed multiscale finite element method using limited global information, simply denoted by MMsFEM, and the mixed generalized
	multiscale finite element method (MGMsFEM) with residual driven online multiscale basis functions. Both methods are under the framework of mixed multiscale finite element methods, where the pressure equation is solved in the coarse grid with carefully constructed multiscale basis functions for the velocity.
	The multiscale basis functions in both methods include local and  global media information. In terms of MsFEM using limited global information, only one multiscale basis function is utilized in each local neighborhood while multiple basis are used in MGMsFEM. We will test and compare these two methods using the benchmark
	three-dimensional SPE10 model. A range of coarse grid sizes and different combinations
	of basis functions (offline and online) will be considered with CPU time reported for each case. In our numerical experiments, we observe good accuracy by the two above methods.
	Finally, we will discuss and compare the advantages and disadvantages of the two methods in terms of accuracy and computational costs.

	\section{Introduction}
	
	Many practical applications such as nuclear waste storage and
	reservoir simulations require multiphase transport simulations
	in porous media.
	Detailed reservoir models with a wide range of scales are available for more accurate simulations. In this paper, we consider a two-phase flow model, in which a pressure equation is coupled with a specific transport equation. As a main component of the model, we consider the permeability that exhibits high heterogeneities with large contrast, which creates a significant difficulty to resolve by traditional approaches. In particular, due to the multiscale nature of the problem, solving on a fine-scale that is sufficient to capture the heterogeneity results considerable computation, which motivates upscaled or multiscale methods and some other model reduction methods. The idea in upscaled method involves using one effective media parameter in each coarse block so that the concerned  equation can be solved in a reduced model \cite{wu2002analysis,DDDAS_upscale_2004,Arbogast_Boyd_06,durlofsky1991numerical,efendiev2002numerical,efendiev2000modeling}. However, this approach may encounter inaccuracy for missing fine-scale information that has strong effect on the overall flow. Another way is to solve the equation
	in coarse grid with multiscale methods \cite{egw10,efendiev2009multiscale,Arbogast_two_scale_04,chung2015mixed,chen2003mixed,jennylt03,Wheeler_mortar_MS_12, Arbogast_PWY_07, mortar_elliptic,mortar_online,yang2019multiscale,yang2018online}.
	Underlying idea in this method is to incorporate small-scale information in the basis functions and then capture this effect on a larger scale. Parallel computations can be used since the construction of multiscale basis in different local neighborhoods is independent, and this enhances the efficiency.
	
	An important consideration for this coupled pressure-convection-diffusion system is that local conservation of mass is required, which can be hardly guaranteed within the framework of continuous Galerkin method. This requirement has motivated a variety of methods, for example, multiscale finite volume methods \cite{cortinovis2014iterative,jenny2003multi,lunati2004multi}, mixed multiscale finite element methods \cite{aarnes2004use,aarnes2008mixed,chen2003mixed,chung2015mixed}, mortar multiscale methods \cite{arbogast2007multiscale, peszynska2005mortar,peszynska2002mortar}, discontinuous Galerkin (DG) methods \cite{du2018adaptive,kim2013staggered,cockburn2002local}, and postprocessing methods \cite{odsaeter2017postprocessing,bush2013application}.  In this paper, we utilize the mixed multiscale finite element method. In a mixed finite element formulation, one uses a first-order system for pressure and velocity, while the pressure space is spanned by piecewise constant pressure basis functions. The support for each pressure basis function is a single coarse block while for velocity basis function, it vanishes outside a coarse neighborhood, which is composed of two coarse blocks sharing a common coarse edge. Contributed by this construction, local conservation in coarse scale can be attained, which is sufficient in most circumstances. On the other hand, when oscillations appear in the source term, conservation property is not satisfied in the relatively small scale, which motivates the need for some postprocessing procedures. In this paper, some relevant methods are reviewed in Section \ref{post}. Since we only need to postprocess the velocity field when the source term is not a constant, the additional computation is limited.
	Furthermore, one advantage of some postprocessing procedures is that the implementation can be conducted in parallel if the corresponding regions are non-overlapping, which improves the efficiency significantly. 
	
	A typical strategy to simulate the multi-phase flow is the implicit pressure explicit saturation (IMPES) method.  That is, the pressure equation  is solved at
	the current time step with mobility coefficients defined with saturation values from
	the previous time step. The saturation equation is then computed with an explicit Euler scheme using velocity field  from the previous time step.
	The dominant computational cost exists in solving the pressure equation.
	Most multiscale methods for multi-phase simulation are trying to reduce the
	computational cost of solving the pressure equation.
	We note that in \cite{kippe2008comparison}, the authors compared several popular
	multiscale methods for the two-phase flow simulation.
	However, no large-scale three-dimensional test was provided.
	
	Our approach is based on  multiscale finite element method (MsFEM), the pioneering work of which can be found in \cite{babuvska1983generalized,hou1997multiscale}, where it was initially proposed to
	solve the elliptic equation in second order formulation. Furthermore, it was then extended to mixed
	multiscale finite element method (MMsFEM). The key ingredient of the method is the construction of multiscale space which contains the fine-grid information of permeability field. Besides, due to the independent nature of the construction of each local multiscale space, we can implement in parallel to save computational time. Based on traditional MMsFEM, a new MMsFEM using limited global information \cite{aarnes2008mixed,jiang2012some} is introduced. Without ambiguity, we use MMsFEM to denote the new method for simplicity. A main attraction of the method is that one may make use of global information obtained by single-phase flow in the construction of multiscale basis functions, while only local information is used in traditional MMsFEM. In particular, we make some modifications in the boundary conditions in the local problems which are aimed at constructing the local multiscale space. This novel modification is powerful when fine-scale feature of two-phase flow depends on single-phase flow. Moreover, since only one velocity multiscale basis function is utilized in each local neighborhood, the simulation time is relatively shortened compared to the case where multiple basis functions are used, which will discussed in the numerical results. A limitation of this method is that one can not flexibly adjust the dimension of approximation space by changing the number of multiscale basis functions, which motivates the other method mixed Generalized Multiscale method (MGMsFEM).
	
	As is mentioned in previous papers \cite{efendiev2013generalized,Efendiev2011ADD,chung2014adaptive},
	one basis function on each edge is not sufficient especially when there are long channels and non-separable scales in the permeability field. MGMsFEM is designed to construct multiple basis functions on some coarse edges, which provides accuracy in this case and guarantee conservation of mass simultaneously. This method is divided into two steps, constructing snapshot space and multiscale space. First, one may construct a snapshot space as approximate space, where there are basis functions carry important properties of solution. Furthermore, a reduced space called offline space is obtained by some well-designed spectral problems. As mentioned in the paper \cite{CHUNG201669}, a good approximation of the reduced model is only attained when the offline information can effectively represents the real case.  At the same time, when the number of offline basis functions exceeds some level, the error decreases will slow down. In particular, the convergence rate is proportional to $1/\Lambda$, where $\Lambda$ is the smallest eigenvalue we abandon during the construction of multiscale space. Consequently, some online enrichment \cite{chan2016adaptive,chung2015residual,online_mixed,online_cg} is necessary for the sake of more accurate velocity fields, which will be inherited in saturation by the following simulation.
	
	The main purpose of the paper is to compare the MGMsFEM with MMsFEM. As is mentioned above, the MMsFEM is efficient as only one multiscale basis function is utilized in each local region, which results in smaller coarse-grid system and hence shorter simulation time. However, it lacks flexibility compared with MGMsFEM as one can adapt the dimension of the approximation space in the latter one method.
	We will show that MGMsFEM serves as an efficient and flexible way to reduce the relevant errors. Different combinations of offline and online basis functions can be selected to provide most satisfying outcomes. Due to the power of the online basis functions, one may use relatively small number of basis functions to get comparable accuracy with MGMsFEM. Moreover, in context of two-phase flow model, there is no need to solve single-phase problem in advance, which is compulsory in MMsFEM. On the other hand, the construction of snapshot space and solving subsequent local spectral problems in MGMsFEM are not necessary in MMsFEM. As for the similarities of these two methods, once multiscale space is constructed, it can be applied to different cases, e.g. source terms, boundary conditions, mobility on coarse grid, etc., and no further update
	will be conducted. It is computationally cheap since flow equations need to be solved multiple times during in multi-phase transport simulation.
	
	The paper is organized as follows. In Section 2, we introduce the two-phase model. Main ingredients of two concerned methods are discussed respectively in the following two sections. In Section 5, the postprocessing procedures are briefly reviewed in Section 5. Finally, some extensive numerical results are exhibited.

	\section{Model formulation}
	In this paper, we consider the two-phase flow model.
	We can write the flow equations as follows:
	\begin{eqnarray}\label{eq:pressure}
	\left\{\begin{aligned}
	(\lambda(S) k(x))^{-1} \mathbf{v}-\nabla p &=0 \quad \text { in } \quad D \\
	\operatorname{div}( \mathbf{v}) &=f \quad \text { in } \quad D \\
	\mathbf{v}\cdot \mathbf{n} &=g(\mathbf{x}) \quad \text { on } \quad \partial D \label{model}
	\end{aligned}\right.
	\end{eqnarray}
	First equation is obtained applying the Darcy's law and the second states the conservation of mass. $\lambda(S)$ is the total mobility. In particular, $\lambda(S)=\lambda_o(S)+\lambda_w(S)$. And we have $\lambda_i(S)=\frac{k_{ri}(S)}{\mu_{i}}$ , for $i=o,w$. The $k_{rw}(S)$ and $k_{ro}(S)$ are permeabiliy for water and oil while $\mu_{w}(S)$ and $\mu_{o}(S)$ are viscosities for water and oil phases, correspondingly. The dynamics of saturation $S$ affects the flow equation. One can write the dynamics of saturation in the following transport equation.
	\begin{eqnarray}
	\dfrac{\partial S}{\partial t}+ \operatorname{div} (\mathbf{v}\cdot f(S))=q\label{transport}
	\end{eqnarray}
	where $f(S)=\frac{\lambda_w(S)}{\lambda_w(S)+\lambda_o(S)}$ and $q$ is some external force term. \\
	To solve (\ref{transport}), we need to solve (\ref{model}) in advance. Solution algorithm for the two-phase problem is presented in table \ref{two-phase algorithm}. \\
	\begin{table}[!htbp]
		\centering
		\begin{tabular}{c l}
			\hline
			&Two-phase algorithm  \\
			\hline
			\textbf{Input} & $S_{n-1}$ obtained in previous time step\\
			\textbf{Output}& $S_{n}$ \\
			&1. Solving (\ref{model}) to get $p_n$ and $\mathbf{v}_n$\\
			&2. Using $\mathbf{v}_n$ and $S_{n-1}$  in (\ref{transport}) to get $S_{n}$ \\
			\hline
		\end{tabular}
		\label{two-phase algorithm}
	\end{table}
	
	To solve (\ref{transport}), we integrate it with $[t_{n-1},t_n]$ and some volumn  $V_i\subset D$
	\begin{eqnarray}
	\text{meas}(V_i)(S_{z,n}-S_{z,n-1})+\Delta t\int_{\partial_{V_i}}\mathbf{v}\cdot n f(S_{z,n-1})d l=\Delta t\int_{V_i} q_w d x
	\end{eqnarray}
	where we have neglected the error terms and
	\begin{eqnarray}
	S_{z,n}\approx\dfrac{1}{meas(V_i)}\int_{V_i}S(x,t_n)d x
	\end{eqnarray}
	We use $meas(A)=\int_{D} 1_{A}d x$ with $1_A=1$ when $x\in A$ while 0 elsewhere.\\
	To evaluate the term $\int_{\partial_{V_i}}\mathbf{v}\cdot n f(S_{z,n-1})d l$, we use an upwinding scheme. A review of upwinding on a rectangular mesh can be found for example in \cite{thomas2013numerical}. It is imperative that the numerical approximation of $\mathbf{v}$ satisfies the following local conservation property. In particular, it is desirable to have
	\begin{eqnarray}
	\int_{\partial{V_i}}\mathbf{v}\cdot \mathbf{n}d l=\int_{V_i} q d x.
	\end{eqnarray}
	We will consider more details about the mass conservation property in Section \ref{post}.
	\section{Mixed GMsFEM}
\begin{figure}[htbp!]
		\centering
		\includegraphics[scale=0.3]{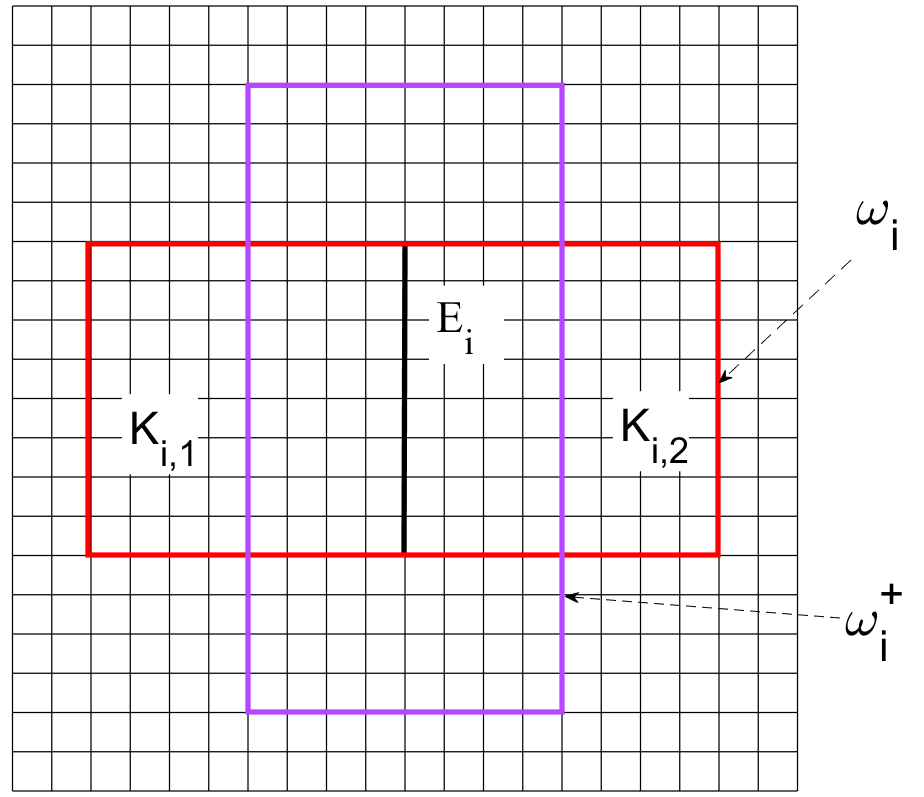}
		\caption{Illustration of coarse neighborhood $\omega_i=K_{i,1}\cup K_{i,1}$ and oversampled neighborhood $\omega_i^+$}
		\label{mixed neighbor}
	\end{figure}

	In this section, we briefly review the Generalized multiscale finite element method for (\ref{eq:pressure}). The construction consists of two steps.
	\begin{enumerate}
		\item Constructing a snapshot space.
		\item Solve a spectral problem.
	\end{enumerate}
	Let $T^{h}$ be a partition of the domain $D$ into fine finite elements. Here $h>0$ is the fine grid mesh size.
	The coarse partition, $T^{H}$ of the domain $D$, is formed such that each element in $T^{H}$ is a connected union of fine-grid blocks. More precisely, $\forall K_{j} \in T^{H}$, $ K_{j}=\bigcup_{F\in I_{j} }F$ for some $I_{j}\subset T^{h}$. The quantity $H>0$ is the coarse-mesh size. We will consider the rectangular coarse elements and the methodology can be used with general coarse elements. An illustration of the mesh notations is shown in the Figure \ref{mixed neighbor}. We denote the interior coarse edges of $T^{H}$ by $E_i,i=1,\cdots,N_{\text{in}}$,
	where $N_\text{in}$ is the number of interior coarse edges. The coarse elements
	of $T^{H}$ are denoted by $K_j,j=1,2,\cdots,N_e$, where $N_e$ is the number of coarse elements. We define the coarse neighborhood of the edge $E_i$ by $\omega_i:=\cup\{K_j\in T_{H}:E_i\subset \overline{K_j}\}$.
	The snapshot space is linearly spanned by a set of extensive basis functions with all possible boundary conditions. We utilize a well-designed local spectral problem to select dominant modes, providing the solution with rapidly decreasing residual, which is illustrated in detail in error analysis section in \cite{efendiev2013generalized}.
	\subsection{Snapshot space}
	We let $V_h$ be the lowest Raviart Thomas vector field ($RT_0$) and $Q_h$ be space of piesewise constant functions on the fine mesh. We set $V_h=\text{span}\{\psi_{\text{fine},1},\cdots,\psi_{\text{fine},N_{\text{in,fine}}}\}$ and
	$Q_h=\text{span}\{p_{\text{fine},1},\cdots,p_{\text{fine},N_{\text{e,fine}}}\}$, where $N_{\text{in,fine}}$ and $N_{\text{e,fine}}$ are the number of inner edges and blocks in fine grid. To construct the multiscale basis functions, we need to construct a snapshot space. In particular, we need to solve the following local problem in each $\omega_i$:
	\begin{eqnarray}
	\left\{
	\begin{aligned}
	(\lambda(S) k(x))^{-1} \mathbf{\psi}^{(i)}_j+\nabla p_j^{(i)}&=0 \quad\text{in}\quad \omega_i\\
	\text{div}(\mathbf{\psi}_j^{(i)})&=\alpha_j^{(i)}\quad \text{in}\quad \omega_i \label{mixed gms}\\
	\mathbf{\psi}_j^{(i)}\cdot n_i&=0\quad \text{on } \partial\omega_i\\
	\mathbf{\psi}_j^{(i)}\cdot m_i&=\delta_j^{(i)}\quad \text{on } E_i \label{snapshot}
	\end{aligned}\right.
	\end{eqnarray}
	$E_i$ is the common edge of the two elements $K_p,K_q\subset\omega_i$. And $E_i=\bigcup_{l_j\subset E_i} l_j$, as is shown in figure \ref{mixed neighbor}.
	\begin{eqnarray}
	\delta_j^{(i)}=
	\left\{
	\begin{aligned}
	&1,\quad\text{on }l_j\\
	&0,\quad\text{elsewhere on }E_i
	\end{aligned}
	\right.
	\end{eqnarray}
	We also set $\alpha_j^{(i)}$ such that the compatibility condition holds:$\int_{K_l} \alpha_j^{(i)}=\int_{E_i}\delta_j^{(i)}$ is satisfied for all $K_l\subset\omega_i$.\\
	In matrix formulation, we solve
	\begin{equation*}
	\left[\begin{array}{cc}
	A_{h}^i& {B_{h}^i}^T \\
	B_{h}^i  & 0
	\end{array}\right]\left[\begin{array}{l}
	\phi_{j,c}^{(i)} \\
	p_{j,c}^{(i)}
	\end{array}\right]=\left[\begin{array}{c}
	0 \\
	\tilde{\alpha}_j^{(i)}
	\end{array}\right],
	\end{equation*}
	where $\phi_{j,c}^{(i)}$ and $p_{j,c}^{(i)})$ are coefficient vectors of $\phi_{j}^{(i)}$ and $p_{j}^{(i)}$ in expansion of space $V_h\times Q_h$. $A_h^i$ and $B_h^i$ are corresponding matrix in $\omega_i$. $\tilde{\alpha}_{i}$ is a vector and in particular,  $\tilde{\alpha}_{i}[j]=\int_{\omega_i} \alpha_i p_i^j$ for all j, where $p_i^j$ is the $j-th$ pressure basis function in space $Q_h$ constrained in $\omega_i$.\\
	We remark that (\ref{mixed gms}) is solved separately in the coarse elements forming $\omega_i$.
	Besides $n_i$ is a unit outer normal vector to $\partial \omega_i$ and $m_i$ is a fixed unit normal vector for $E_i$.\\
	From (\ref{snapshot}), we can solve snapshot basis $\mathbf{\psi}_j^{(i)}$ for $1\leq j\leq J_i$.
	We can further define $V_{\text{snap}}^{i}=\text{span}
	\{\mathbf{\psi}_1^{(i)},\cdots,\mathbf{\psi}_{J(i)}^{(i)}\}$ for $1\leq i\leq N_{\text{in}}$, where $N_{\text{in}}$ is the number of inner edges for coarse elements. Finally, $V_{\text{snap}}=\bigoplus_{i=1}^{N_{\text{in}}}V_{\text{snap}}^{i}$.
	
	\subsection{Offline space}
	We proceed to the space reduction in this subsection. We aim at solving the eigenpair $(\lambda,\mathbf{v})$ from the following spectral problem,
	\begin{eqnarray}
	a(\mathbf{v},\mathbf{w})=\lambda s(\mathbf{v},\mathbf{w}),\label{spectral}
	\end{eqnarray}
	where $a(\mathbf{v},\mathbf{w})$ and $s(\mathbf{v},\mathbf{w})$ are symmetric positive definite bilinear forms defined on $V_{\text{snap}}\times V_{\text{snap}}$. We can view the spectral problem as a way to select the dominating components of a specific residual operator. In particular, we can define $\mathcal{A}:V_{\text{snap}}\rightarrow V_{\text{snap}}$ as follows.
	\begin{eqnarray}
	s(\mathcal{A}\mathbf{v},\mathbf{w})=a(\mathbf{v},\mathbf{w})
	\end{eqnarray}
	where $s(\mathbf{v},\mathbf{w})$ is an inner product on $V_{\text{snap}}$. We assume $\mathcal{A}$ has a rapidly decaying eigenvalues. In practice, we need to solve a local spectral problem, which is more sufficient than the global one. In particular, in $\omega_i$ corresponding to $E_i$ we find a real number $\lambda$ via solving the following problem:
	\begin{eqnarray}
	a_i(\mathbf{v},\mathbf{w})=\lambda s_i(\mathbf{v},\mathbf{w})\label{spectral_i}
	\end{eqnarray}
	where $a_i(\mathbf{v},\mathbf{w})=\int_{E_i}(\lambda(S) k(x))^{-1}(\mathbf{v}\cdot m_i)(\mathbf{w}\cdot m_i)$ and $s_i(\mathbf{v},\mathbf{w})=\int_{\omega_i}(\lambda(S) k(x))^{-1}\mathbf{v}\cdot \mathbf{w}+\int_{\omega_i}\text{div}(\mathbf{v})\cdot\text{div}(\mathbf{w})$. \\
	Suppose we have solved $\lambda_{1}^{(i)}\cdots\lambda_{J_i}^{(i)}$ corresponding to eigenvectors $\mathbf{v}_1^{(i)}\cdots \mathbf{v}_{J(i)}^{(i)}$ in $\omega_i$, we further select first $l_i$ eigenvectors and consequently the offline basis function
	$\mathbf{\phi}_j^{(i)}=\sum_{k=1}^{J_i} v_{jk}^{(i)}\psi_k^{i}$, where $v_{jk}^{(i)}$
	is the k-th component of $\mathbf{v}_{j}^{(i)}$.\\
	And we can construct local offline space in $\omega_i$.
	$V_{\text{off}}^{(i)}=\text{span}\{\mathbf{\phi}_1^{(i)},\cdots,\mathbf{\phi}_{l_i}^{(i)}\}$.
	We set the global offline basis function to be $V_{\text{off}}=\bigoplus_{E_{l} \in \mathcal{E}_{n}} V_{\text{off}}^{i}$, where $\mathcal{E}_{n}$ is the collection of all inner coarse faces.
	\subsection{Online enrichment}
	To improve the accuracy of velocity solution in (\ref{model}), we consider enriching the offline space based on some residual functions. Our enrichment is performed iteratively. Based on the space $V_{H}^l$ from previous iteration, we can obtain $V_{H}^{l+1}$ through the following 5 steps. In particular, we have $V_{H}^0=V_{\text{off}}$.
	First, we define a residual operator on this iteration $$R_{D}^l(\mathbf{v})=\int_{D} (\lambda(S) k(x))^{-1}\mathbf{v}_H^l \cdot \mathbf{v}-\int_{D} \text{div}(\mathbf{v})p_H^l\quad \forall
	v\in V_{\text{snap}},$$
	
	Let $\hat{V}_{\text{snap}}$ be the divergence free subspace of ${V}_{\text{snap}}$, and we have
	$$R_{D}(\mathbf{v})=\int_{D} (\lambda(S) k(x))^{-1}\mathbf{v}_H \cdot \mathbf{v}\quad \forall
	\mathbf{v}\in\hat{V}_{\text{snap}}.$$\\
	Step 1: Find the multiscale solution  $\mathbf{v}_{H}^{l} \in V_{H}^{l}$ and $p_{H}^{l} \in Q_{H}$ satisfying
	\[
	\begin{aligned}
	\int_{D} (\lambda(S) k(x))^{-1} \mathbf{v}_{H}^{l} \cdot \mathbf{w}-\int_{D} \operatorname{div}(\mathbf{w}) p_{H}^{l} &=0 \quad \forall w \in V_{H}^{l}, \\
	\int_{D} \operatorname{div}\left(\mathbf{v}_{H}^{l}\right) q &=\int_{D} f q \quad \forall q \in Q_{H}.
	\end{aligned}
	\]
	In matrix formulation, we have
	\begin{equation}
	\left[\begin{array}{cc}
	\left(V^{l}\right)^{T} A_{h} V^{l} & \left(V^{l}\right)^{T} B_{h}^{T} P_{H} \\
	P_{H}^{T} B_{h} V^{l} & 0
	\end{array}\right]\left[\begin{array}{l}
	\mathbf{v}_{H, c}^{l} \\
	\mathbf{p}_{H, c}^{l}
	\end{array}\right]=\left[\begin{array}{c}
	0 \\
	P_{H}^{T} F_{h}
	\end{array}\right].\label{matrix gms}
	\end{equation}
	Step 2:  Let $\{\omega_i\}$ be a non-overlapping partition of $D$, which corresponds to common coarse faces $\{E_i\}$. We let $\{\omega_i^+\}$ be oversampled neighborhoods from $\{\omega_i\}$. We compute $R_{\omega_{i}^{+}}(\mathbf{v})$ for each $\mathbf{v}\in V_{\text {snap }}^{i ,+}$. We note that we need to project the $\mathbf{v}_{H}^{l}$ and $p_{H}^{l}$ into the fine-grid spaces $V_{h}$ and $Q_{h}$, respectively, in order to compute $R_{\omega_{i}^{+}}^l(v)$.  \\
	Step 3: For each $\omega_{i}^{+},$ we solve for $\mathbf{\phi}_{i}^{+} \in V_{\text {snap }}^{i ,+}$ such that
	\[
	\int_{\omega_{i}^{+}} (\lambda(S) k(x))^{-1} \mathbf{\phi}_{i}^{+} \cdot \mathbf{v}=R_{\omega_{i}^{l,+}}(\mathbf{v}) \quad \forall \mathbf{v} \in V_{\text {snap }}^{i ,+}
	\]
	\\
	Step 4: We take the restriction of $\mathbf{\phi}_{i}^{+} \cdot n_{i}^{+}$ on the coarse face $E_{i}$ and normalize it, and we denote it by $\lambda_{i}$.\\
	Step 5: The online basis $\mathbf{\chi}_{i}$ for the coarse neighborhood $\omega_{i}$ can be constructed by solving
	\begin{eqnarray}
	\begin{aligned}
	(\lambda(S) k(x))^{-1} \mathbf{\chi}_{i}+\nabla p_{i} &=0 \quad \text { in } \omega_{i} \\
	\operatorname{div}\left(\mathbf{\chi}_{i}\right) &=\alpha_{i} \quad \text { in } \omega_{i} \\
	\mathbf{\chi}_{i} \cdot n_{i} &=\lambda_{i} \quad \text { on } E_{i} \\
	\mathbf{\chi}_{i} \cdot n_{i} &=0 \quad \text { on } \partial \omega_{i} \label{online glo}
	\end{aligned}
	\end{eqnarray}
	where $\alpha_{i}$ is chosen to satisfy the condition $\int_{K} \alpha_{i}=\int_{E_i} \lambda_i$ for every $K \subseteq \omega_{i}, n_{i}$ is a fixed unit-normal vector for the coarse face $E_{i} .$\\
	In matrix formulation, we have
	\begin{equation*}
	\left[\begin{array}{cc}
	A_{h}^i& {B_{h}^i}^T \\
	B_{h}^i  & 0
	\end{array}\right]\left[\begin{array}{l}
	\chi_{i,c} \\
	p_{i,c}
	\end{array}\right]=\left[\begin{array}{c}
	0 \\
	\tilde{\alpha}_{i}
	\end{array}\right],
	\end{equation*}
	where $\chi_{i,c}$ and $p_{i,c}$ are coefficient vectors of $\chi_{i}$ and $p_{i}$ in expansion of $V_h\times Q_h$ constrained in $\omega_i$.\\
	Those $\{\mathbf{\chi}_i\}_{i=1}^{N_{\text{in}}}$ are the new online basis functions. We update the velocity space by letting $V_{H}^{l+1}=V_{H}^l \cup \operatorname{span}\{\mathbf{\chi}_{1},\mathbf{\chi}_{2}, \ldots, \mathbf{\chi}_{N_{\text{in}}}\}$.\\

	\section{Multiscale Finite element Method with limited global information}
	Different from standard Multiscale Finite Element Method, we use construct the multiscale basis function using the single-phase flow information. In particular, we can split this procedure into the following two steps.
	\begin{enumerate}
		\item Step1: Using standard finite element method to solve the fine-scale solution $\mathbf{v}_{sp}$ for the single-phase flow system.
		\item Step2: Construct the Multiscale basis using the global limited information $\mathbf{v}_{sp}$.
	\end{enumerate}
	
	For convenience, we first rewrite the two-phase flow equation as:
	\begin{eqnarray}
	\left\{\begin{aligned}
	(\lambda(S) k(x))^{-1} \mathbf{v}-\nabla p &=0 \quad \text { in } \quad D, \\
	\operatorname{div}( \mathbf{v}) &=f \quad \text { in } \quad D, \\
	\mathbf{v}\cdot \mathbf{n} &=g(x) \quad \text { on } \quad \partial D.
	\end{aligned}\right.
	\end{eqnarray}
	More specifically, the system is reduced to single-phase when $\lambda(S)=1$.
	
	The variational formulation associated with (\ref{model}) is to seek $(\mathbf{v},p)\in H(\text{div},\Omega)\times L^2(\Omega)$ such that $\mathbf{v}\cdot\mathbf{n}=g$ on $\partial \Omega$ and
	\begin{eqnarray}
	\begin{aligned}
	\int_{D} (\lambda(S) k(x))^{-1} \mathbf{v}\cdot \mathbf{w}-\int_{D} \text{div}(\mathbf{w})p_h&=0 \quad\forall \mathbf{w}\in \mathbf{V}_{0},\\
	\int_D \text{div}(\mathbf{v})q&=\int_D f q \quad\forall q\in Q.\label{varia}
	\end{aligned}
	\end{eqnarray}
	where
	$$\mathbf{V}_{0}=\{\mathbf{v}\in H(\text{div},\Omega)|\mathbf{v}\cdot\mathbf{n}=0\}\text{ and } Q= L^2(\Omega).$$
	
	We denote
	$$a(\mathbf{u},\mathbf{v})=((\lambda k)^{-1}\mathbf{u},\mathbf{v}),\quad b(\mathbf{v},p)=(\text{div}(\mathbf{v}),p).$$
	
	Then we can rewrite the above system as below:
	\begin{eqnarray}
	\begin{aligned}
	a(\mathbf{v},\mathbf{w})+b(\mathbf{w},p)&=0 \quad\forall \mathbf{w}\in \mathbf{V}_{0},\\
	b(\mathbf{v},q)&=(f,q)\quad\forall q\in Q.\label{varia matrix}
	\end{aligned}
	\end{eqnarray}
	
	We denote the $(\mathbf{v}_{h},p_h)$ to be the reference solution in $(V_h,Q_h)$ by solving
	\begin{eqnarray}
	\begin{aligned}
	a(\mathbf{v}_{h},\mathbf{w})+b(\mathbf{w},p_{h})&=0 \quad\forall \mathbf{w}\in \mathbf{V}_h^0,\\
	b(\mathbf{v}_{h},q)&=(f,q)\quad\forall q\in Q_h.\label{fine matrix}
	\end{aligned}
	\end{eqnarray}
	where $\mathbf{v}_{h}\cdot n=g_h$ on $\partial D$ and $\mathbf{V}_h^0=V_h\cap\{\mathbf{v}\in V_h:\mathbf{v}\cdot n=0 \text{ on } \partial D\}$.
	
	We can also write (\ref{fine matrix}) in the following matrix form:
	\begin{equation*}
	\left[\begin{array}{cc}
	A_{h}& B_{h}^T \\
	B_{h}  & 0
	\end{array}\right]\left[\begin{array}{l}
	\mathbf{v}_{h,c} \\
	\mathbf{p}_{h,c}
	\end{array}\right]=\left[\begin{array}{c}
	0 \\
	F_{h}
	\end{array}\right]
	\end{equation*}

	where we let $v_{h,c}$ and $p_{h,c}$ to be the corresponding coefficient vector for reference solution $v_h$ and $p_h$. $A_h[i,j]=a(\psi_{\text{fine},i},\psi_{\text{fine},j})$ and $B_h[i,j]=b(\psi_{\text{fine},i},p_{\text{fine},j})$.
	
	We first seek $\mathbf{v}_{sp}$ for the single-phase system in  $\mathbf{V}_h$ by solving (\ref{fine matrix}) with $\lambda(S)=1$.
	
	Then we consider the constructing multiscale basis functions by using the $v_h$ obtained above. We compute $\phi_{i}^K$ by solving the following local problem on coarse neighborhood $\omega_i$.
	\begin{eqnarray}
	\left\{\begin{aligned}
	(\lambda(S) k(x))^{-1} \mathbf{\phi}_{i}^K+\nabla(p_{i}^K)&=0 \text{ in } \omega_i\\
	\text{div}(\mathbf{\phi}_{i}^K)&=\alpha_i \text{ in } \omega_i\\
	\mathbf{\phi}_{i}^K\cdot \mathbf{n}_i&= \mathbf{v}_{sp}\cdot n_i\text{ on } E_i\\
	\mathbf{\phi}_{i}^K\cdot \mathbf{n}_i&=0 \text{ on } \partial{\omega_i}
	\end{aligned}\right.
	\end{eqnarray}
	Remark that $\alpha_{i}$ is chosen to satisfy the condition $\int_{K} \alpha_{i}=\int_{E_i} \lambda_i$ for every $K \subseteq \omega_{i}.$ And the multiscale finite element space $\mathbf{V}_{ms}$ for velocity will be defined by
	\begin{eqnarray}
	\begin{aligned}
	\mathbf{V}_{ms} &=\bigoplus_{K}\left\{\mathbf{\phi}_{i}^K\right\} \subset H(d i v, \Omega), \\
	\mathbf{V}_{ms}^{0} &=\mathbf{V}_{ms} \cap H_{0}(d i v, \Omega).
	\end{aligned}
	\end{eqnarray}
	We then seek the multiscale velocity solution in $(\mathbf{v}_{\text{ms}},p_{\text{ms}})\in (\mathbf{V}_{\text{ms}},Q_{\text{ms}})$ by solving the following
	\begin{eqnarray}
	\begin{aligned}
	a(\mathbf{v}_{\text{ms}},\mathbf{w})+b(\mathbf{w},p_{\text{ms}})\quad\forall \mathbf{w}\in \mathbf{V}_{\text{ms}},\\
	b(\mathbf{v}_{\text{ms}},q)=(f,q)\quad\forall q\in Q_{\text{ms}},\label{coar matrix}
	\end{aligned}
	\end{eqnarray}
	where $Q_{\text{ms}}$ is spanned by piecewise constant functions.
	
	Again, we can use matrix representation in the following.
	\begin{equation}
	\left[\begin{array}{cc}
	\left(V_{Hc}\right)^{T} A_{h} V_{Hc} & \left(V_{Hc}\right)^{T} B_{h}^{T} P_{Hc} \\
	P_{Hc}^{T} B_{h} V_{Hc} & 0
	\end{array}\right]\left[\begin{array}{l}
	\mathbf{v}_{H,c}^{l} \\
	p_{H,c}^{l}
	\end{array}\right]=\left[\begin{array}{c}
	0 \\
	P_{H}^{T} F_{h}
	\end{array}\right].\label{matrix lim}
	\end{equation}
	where $V_{H,c}$ and $P_{H,c}$ are the matrix form of the velocity multiscale space and pressure multiscale space. In particular, we have $V_{H,c}=\text{span}\{\phi_1,\cdots,\phi_n\}$ and $P_{H,c}=\text{span}\{p_1,\cdots,p_m\}$, where $n$ and $m$ are the dimension for velocity multiscale space and pressure multiscale space.
	\subsection{Some remarks about MsFEM with limited global information}
	As in \cite{jiang2012some}, the $(\mathbf{v}_{\text{ms}},p_{\text{ms}})$ is exact for the single-phase system. For the error analysis, we first need some assumptions as below:
	\begin{eqnarray}
	a(\mathbf{u}_{\text{ms}},\mathbf{u}_{\text{ms}})\text{ is  ker} B_{\text{ms}} -coercive \label{assump1}
	\end{eqnarray}
	\begin{eqnarray}
	\inf _{\mathrm{q} \in Q_{\text{ms}}} \sup _{\mathbf{v}\in \mathrm{V}_{\text{ms}}} \frac{b\left(\mathbf{v}, q\right)}{\left\|\mathbf{v}\right\|_{H(d i v, \Omega)}\left\|q\right\|_{L^{2}(\Omega)}} \geq C.
	\label{assump2}
	\end{eqnarray}
	Then we will have the following approximation property:
	\begin{lemma}
		If $(\mathbf{v}_h,p)$ and $(\mathbf{v}_{\text{ms}},p_{\text{ms}})$ are solutions for (\ref{fine matrix}) and (\ref{coar matrix}) respectively and the conditions (\ref{assump1}), (\ref{assump2}) hold, then we have
		\begin{eqnarray}
		\left\|\mathbf{v}_h-\mathbf{v}_{\text{ms}}\right\|_{H(\text {div}, \Omega)}+\left\|p_h-p_{\text{ms}}\right\|_{0, \Omega} \leq \displaystyle\inf_{
			\begin{array}{cc}
			&\mathbf{u} \in \mathbf{V}_{\text{ms}}\\
			&\mathbf{u}-g_{\text{ms}} \in \mathbf{V}^0_{\text{ms}}
			\end{array}
		}\left\|\mathbf{v}_h-\mathbf{u}\right\|_{H(\text {div}, \Omega)}+\inf _{q\in Q_{\text{ms}}}\left\|p_h-q\right\|_{0, \Omega},
		\end{eqnarray}
		where $g_{\text{ms}}$ is the multiscale interpolation in the space $\mathbf{V}_{\text{ms}}$.
	\end{lemma}
	Contributed by this property, we aim at seeking a good enough velocity solution in approximation space to have satisfying approximation effect.
	$(\mathbf{v}_h,p_h)$ and $(\mathbf{v}_{\text{ms}},p_{\text{ms}})$ are solutions for (\ref{fine matrix}) and (\ref{coar matrix}) respectively, we have
	\begin{theorem}\label{the}
		\begin{eqnarray}
		\left\|\mathbf{v}_h-\mathbf{v}_{\text{ms}}\right\|_{H(\text {div}, \Omega)}+\left\|p_h-p_{\text{ms}}\right\|_{0, \Omega}\leq C\delta+CH,
		\label{thm}
		\end{eqnarray}
		where $\delta$ is a very small number and $H$ is the size of the coarse grid.
	\end{theorem}
	Consequently, we can base our multiscale space for two-phase model on the single-phase reference solution to get relatively satisfying estimation.
	\section{Postprocessing procedures}\label{post}
	In mixed formulation, contributed by fact that $Q_H$ is spanned by piecewise constant basis functions, we can have mass conservation in coarse scale as shown in (\ref{coarse conser}), which is not sufficient under some specific circumstances, for example, when the source term has some oscillations in fine scale. In order to obtain the conservation in fine scale, we can make use of postprocessing procedures introduced in \cite{chung2015mixed,guiraldello2020velocity}. \\
	\begin{eqnarray}
	\int_{K} \operatorname{div}\left(\mathbf{v}_{H}\right)  &=\int_{K} f  \quad \forall K\subset D.\label{coarse conser}
	\end{eqnarray}
	In \cite{guiraldello2020velocity}, three types of approaches to postprocess the flux are mentioned, the Mean method, the Patch method and the Stitch method. The latter two are developed from the first one, which we will give a brief review as follows. The basic idea is that after we obtained the multiscale solution $\mathbf{v}_H$, we further solve a new velocity field $\mathbf{v}_{h,p}$ by solving a downscaled local problem, where previous solution information is contained in boundary condition as shown in (\ref{fine conser}).
	\begin{eqnarray}
	\begin{aligned}
	\int_{K} (\lambda(S) k(x))^{-1}\mathbf{v}_{h,p}\cdot w_h-\int_{D} \text{div}(\mathbf{w}_h)p_h&=0 \quad\forall w_h\in V_h(K),\\
	\int_K \text{div}(\mathbf{v}_{h,p})q_h&=\int_D f q_h \quad\forall q_h\in Q_h(K).\label{fine conser}\\
	\mathbf{v}_{h,p}\cdot \mathbf{n} &=\mathbf{v}_{H}\cdot \mathbf{n} \quad \partial K.
	\end{aligned}
	\end{eqnarray}
	Since different coarse blocks are disjoint,implementations can be conducted in parallel to enhance efficiency. In the Patch method, there is a modification in the boundary condition, where we need to solve an auxiliary field in some local problems with small size, but in return higher accuracy is obtained with similar complexity. In terms of Stitch method, on the other hand, computational cost is reduced by restricting the local problems to smaller size.
	
	In some flow problems such as single-phase flow and multi-phase flow model, we will utility the postprocessing techniques to guarantee the mass conservation in fine scale. Despite the fact that computation in fine scale is demanding, this approach is necessary when the source term is non-constant, which means the computations of $v_{h,p}$ are very limited and efficient.
	\section{Numerical results}
\begin{figure}
		\centering
		\includegraphics[scale=0.5]{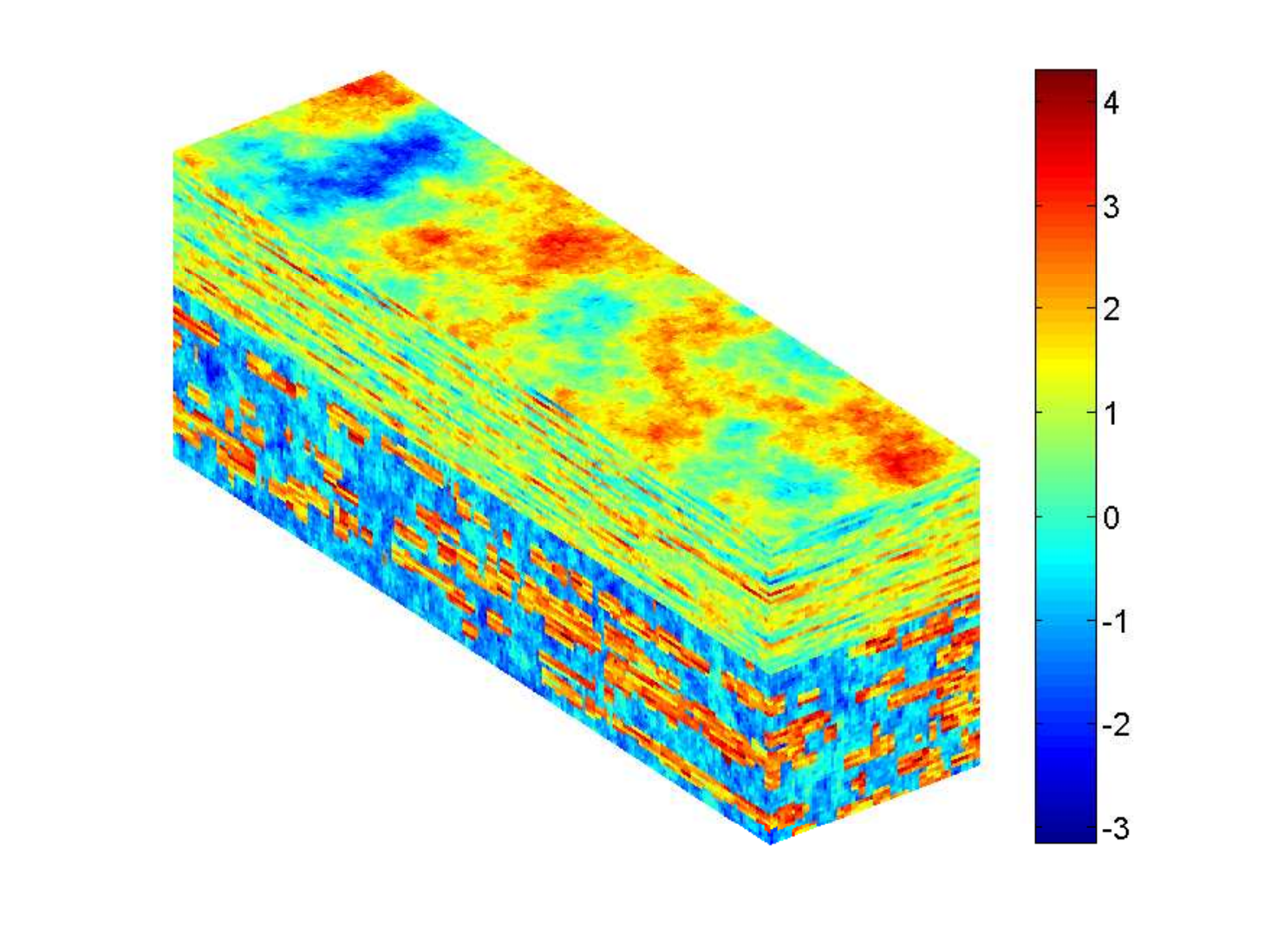}
		\caption{SPE10 model in log10 scale.}
		\label{fig:model}
	\end{figure}

	In this section, we present numerical results for two-phase flow problem, where permeability field comes from SPE10. As one may refer Figure \ref{fig:model}, the most significant features of this permeability are the channelized structure and high contrast, which give rise to some inaccuracy using local approaches. In our experiments, we take the last 80 layers, which contain most heterogeneous part of information. One can observe that it performs better to utilize global information in constructing multiscale basis. In our experiments, We compare the results of Mixed finite element method (MFEM), Mixed multiscale finite element method with limited global information (MMsFEM), Mixed generalized multiscale finite element method (MGMsFEM).
	We performed all the computation on a workstation with Intel(R) Xeon(R) CPU E5-2650 CPU and Matlab, where 40 cores were used in the basis computation stage in parallel.
	
	The permeability field $\kappa(x)$ is given on fine mesh with size $220\times60\times80$ while different scales of coarse mesh are used. In the simulation, we consider $n=20,10,5$, where $n$ is the number of fine elements in each coarse element in each axis. In particular, there are $n^3$ fine elements in each coarse element. We compare two circumstances, where the difference lies in location of
	injectors and producers. In case 1 the water is injected from four edges and produced in the center of the three-dimensional domain while in case 2 we exchange the injectors ans producers. We denote four producers by producer 1, producer 2, producer 3 and producer 4.
	
	We use $L^2$ error for the average relative saturation error $e_s$, in particular, we have
	\begin{eqnarray*}
		e_{s,i}&=\sqrt{\dfrac{\int_{D}(s_{ref,i}-s_{ms,i})^2}{\int_{D}{s_{ref,i}}^2}},\quad     e_s=\dfrac{\sum_{i}^{T_{\text{final}}}e_{s,i}}{n},
	\end{eqnarray*}
	where $s_{ref,i}$ and $s_{ms,i}$ are reference and multiscale saturation at $i-th$ time instant respectively while $T_{\text{final}}$ is the final time instant.
	
	In Table \ref{time1} and Table \ref{time2} (for case 1 and 2 respectively), time in seconds for constructing multiscale basis ($T_{\text{setup}}$ in the table), simulation time and average saturation error $e_s$ among iterations are exhibited. The shown $T_{\text{setup}}$ in two tables can be shortened if more cores are used.
	"Dof" corresponds to the dimension of (\ref{matrix gms}) or (\ref{matrix lim}). Specifically, in each coarse element, we use one pressure basis function which is a piecewise constant function support in this coarse element. As for velocity basis functions, we construct them in each coarse neighborhood. For MMsFEM, there is only one for each coarse neighborhood while multiple basis functions are constructed in MGMsFEM. We use $a+b$ to denote there are $a$ offline and $b$ online basis functions in each local neighborhood. We compare six different combinations of basis functions, $1+0,6+0,8+0,2+1,2+2,4+2$. With same $n=20$, there is a significant improvement in accuracy when we enrich the multiscale space with online basis functions. The average errors shown in the lase column of Table \ref{time1} are approximately diminished by half from $8\%$ in the case $8+0$ to $4\%$ in case $4+2$. In the second example, there is a degradation in accuracy after injectors and producers are exchanged. However, the error decay is still apparent especially when online basis functions are incorporated. The case $2+2$ has smaller average error $7.7\%$ compared with $8.6\%$ in the case $8+0$, where the corresponding dimension is nearly doubled and simulation time also increased by almost a half, nearly 1500 seconds. The case $4+2$ with $n=20$ performed practically as well as the case $2+2$ with $n=10$, while in the latter case the dimension is much larger, which gives rise to almost twice the length of simulation time, i.e. more than 8000 seconds. Consequently, it is evident that the enrichment, in both offline and online levels, can improve the approximation performance. However, online enrichment is more efficient than offline enrichment, where only local information is utilized in the latter case.
	
\begin{figure}[H]
		\centering	
		\includegraphics[trim={1cm 6cm 1cm 7cm},clip,width=3.8in]{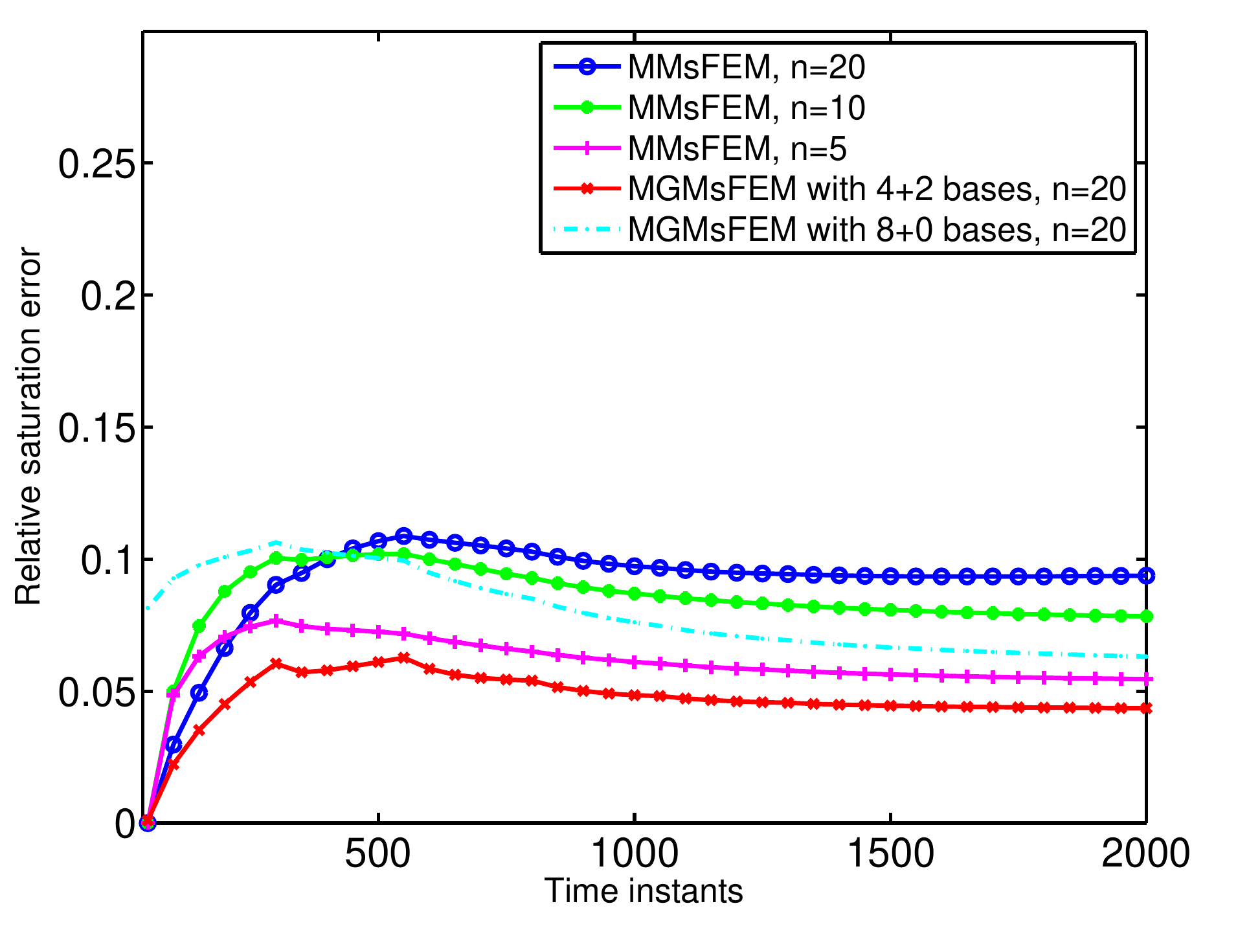}
		\caption{Relative saturation errors as functions of time for case 1.}
		\label{fig:e_1}
	\end{figure}

\begin{figure}[H]
		\centering	
		\includegraphics[trim={1cm 6cm 1cm 7cm},clip,width=3.8in]{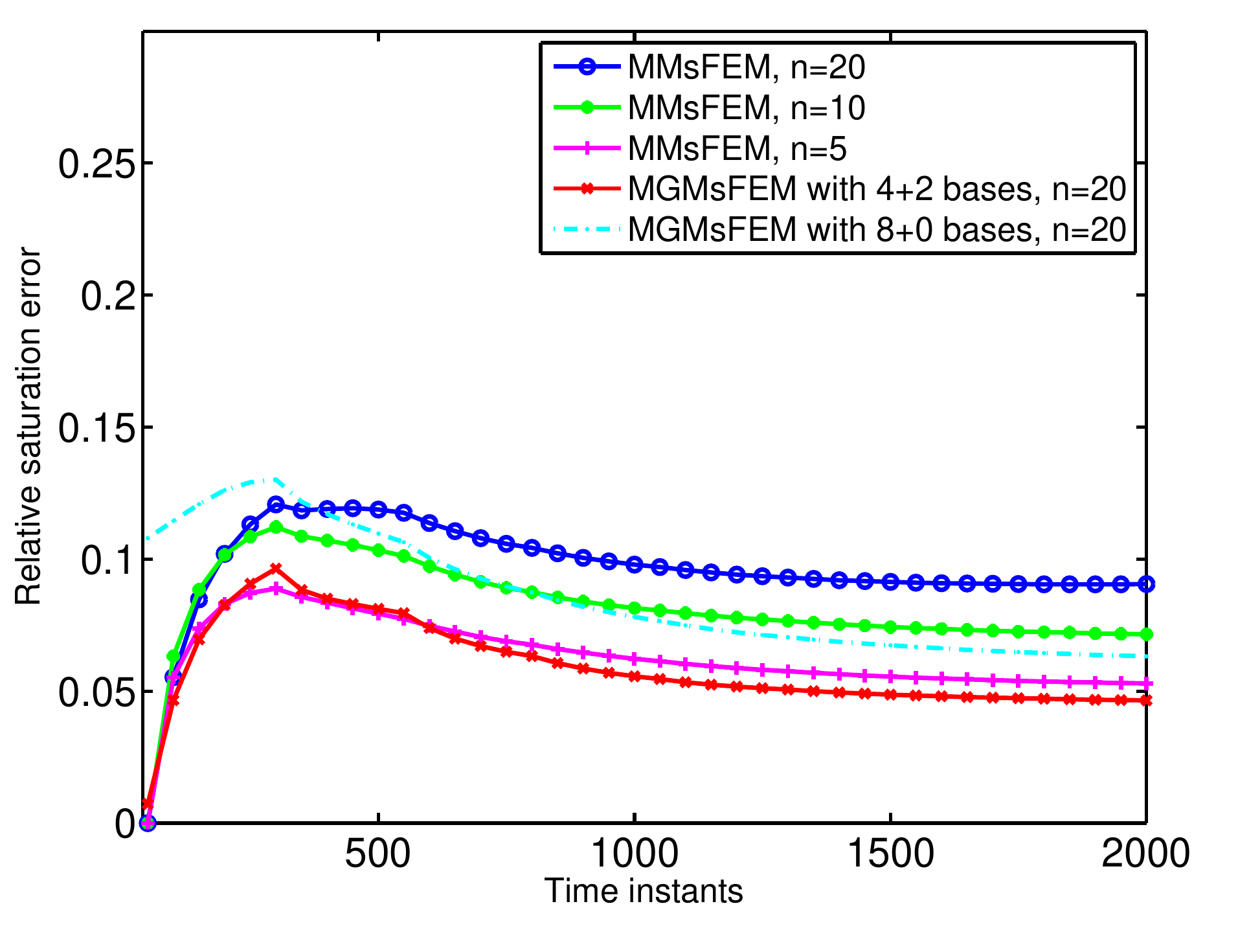}
		\caption{Relative saturation errors as functions of time for case 2}
		\label{fig:e_2}
	\end{figure}

	As for error convergence, one may easily verify from Table \ref{time1} and Table \ref{time2} that error decay is apparent with finer coarse elements for MMsFEM. We consider three sizes of coarse elements, $n=20,10,5$. Notably, when $n$ is twice larger, the simulation time is much shorter without losing much accuracy. For instance, the average error is increased by $2\%$ when we resize the coarse elements from $n=5$ to $n=10$ while the dimension is largely diminished. In addition, the simulation time is reduced by nearly $80\%$, which contributes to higher efficiency. In Theorem \ref{the}, there is error bound $\delta$, which represents how well the two-phase velocity field can be approximated by the single-phase velocity field in each local patch. As we decrease the coarse-mesh size, the error will converge to the bound. From Figure \ref{fig:e_1} and \ref{fig:e_2}, it is remarkable that the largest saturation error is nearly reduced by half in $4+2$ compared to $8+0$, which is consistent with the aforementioned average error. In addition, one can observe that it converges faster in iterations with bigger dimension. In particular, before 300 time instants, the accuracy in case $n=5$ is almost steady while for case $n=20$, the error keeps increasing in the first 500 instants.
	
	When we compare the MMsFEM and MGMsFEM, we can see that both of these two methods can provides relatively accurate solutions. Compared to the case when only offline basis functions are used in MGMsFEM, MMsFEM is more efficient in reducing error while the dimension of the system is much smaller than MGMsFEM, which costs shorter computation time. It is apparent that limited global information efficiently improves the accuracy without increasing much computation time by comparing the case between MMsFEM and MGMsFEM with $1+0$ basis. For example when $n=20$, MMsFEM gives $9.1\%$ compared with $39.3\%$ MGMsFEM with $1+0$ basis while the time spent in multiscale space construction and simulation is similar in these two cases. In addition, when we consider the case when $n=20$, the MMsFEM has average error $9\%$ compared to $8\%$ in MGMsFEM(8+0). However, the simulation time in MMsFEM is shortened by half (about 3000 seconds). Once the online enrichment is implemented in MGMsFEM, error decay is remarkable while computation time is comparable with MMsFEM. In Table \ref{time1}, when $n=20$, we can see that MMsFEM has error $9.1\%$ with near 3363 seconds while MGMsFEM$(2+1)$ has error $11.2\%$ with about 3700 seconds. However, MGMsFEM$(2+2)$ achieves smaller error $6.4\%$ with not much longer simulation time about 4000 seconds. Consequently, we can come to the conclusion that MGMsFEM serves as a flexible tool as one can adjust the numbers of offline and online basis functions to achieve relatively high accuracy without adding some computation costs.
	
\begin{figure}[H]
		\centering
		\subfigure[reference]{
			\includegraphics[trim={1cm 7.5cm 1cm 7cm},clip,width=3.0in]{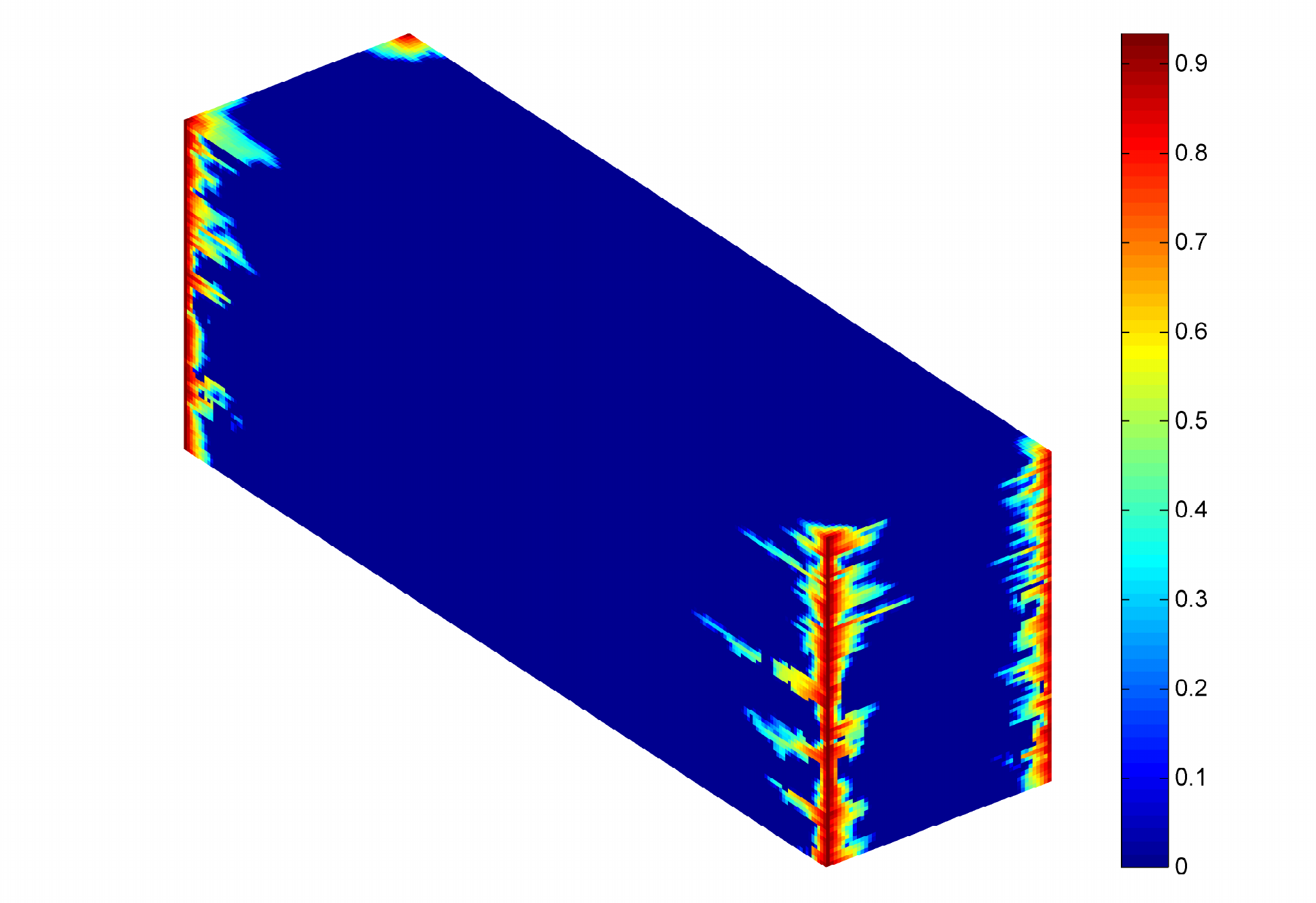}}
		\subfigure[MMsFEM with n=10]{
			\includegraphics[trim={1cm 7.5cm 1cm 7cm},clip,width=3.0in]{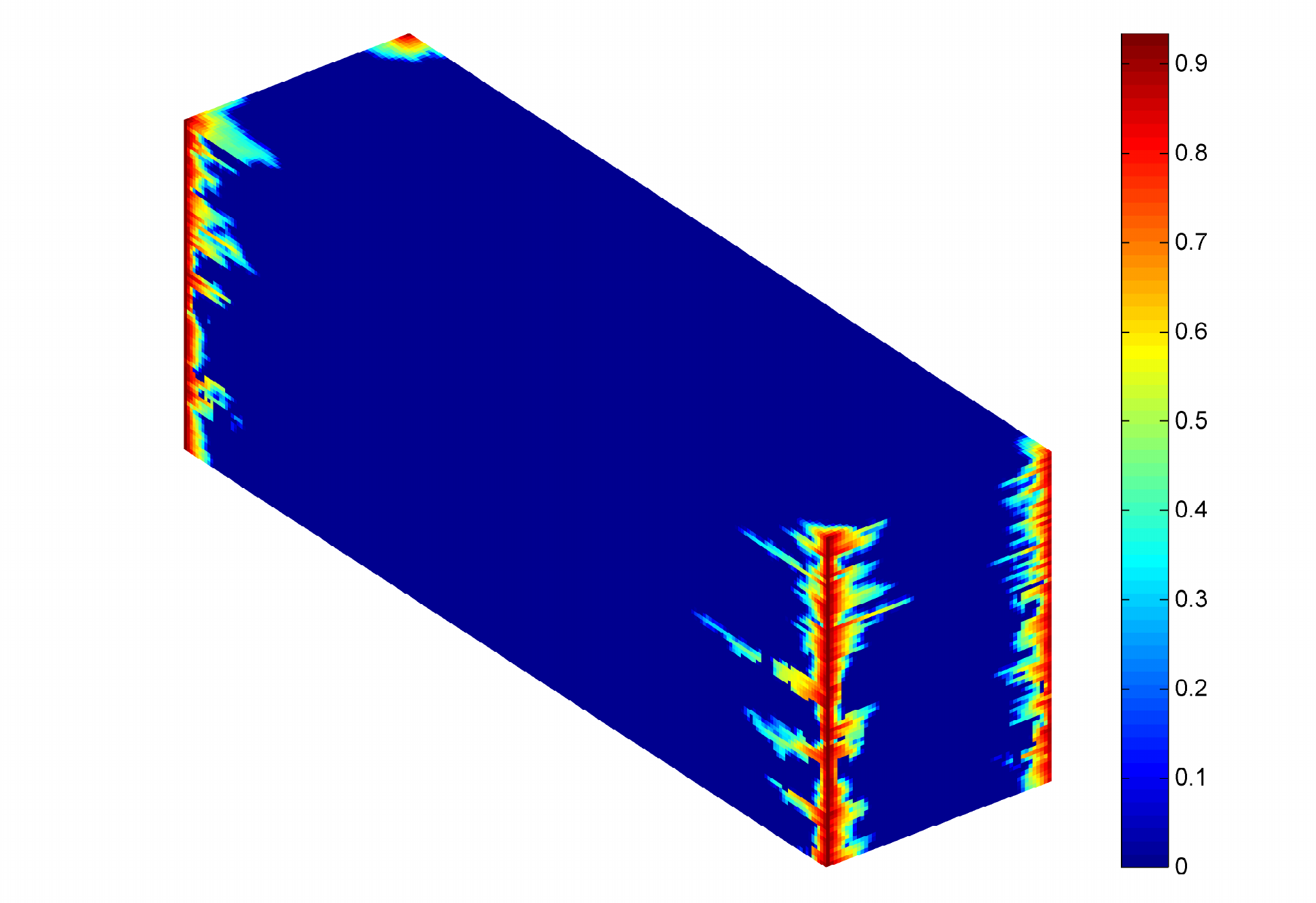}}
		\subfigure[MGMsFEM with n=20 and 4+2 basis functions]{
			\includegraphics[trim={1cm 7.5cm 1cm 7cm},clip,width=3.0in]{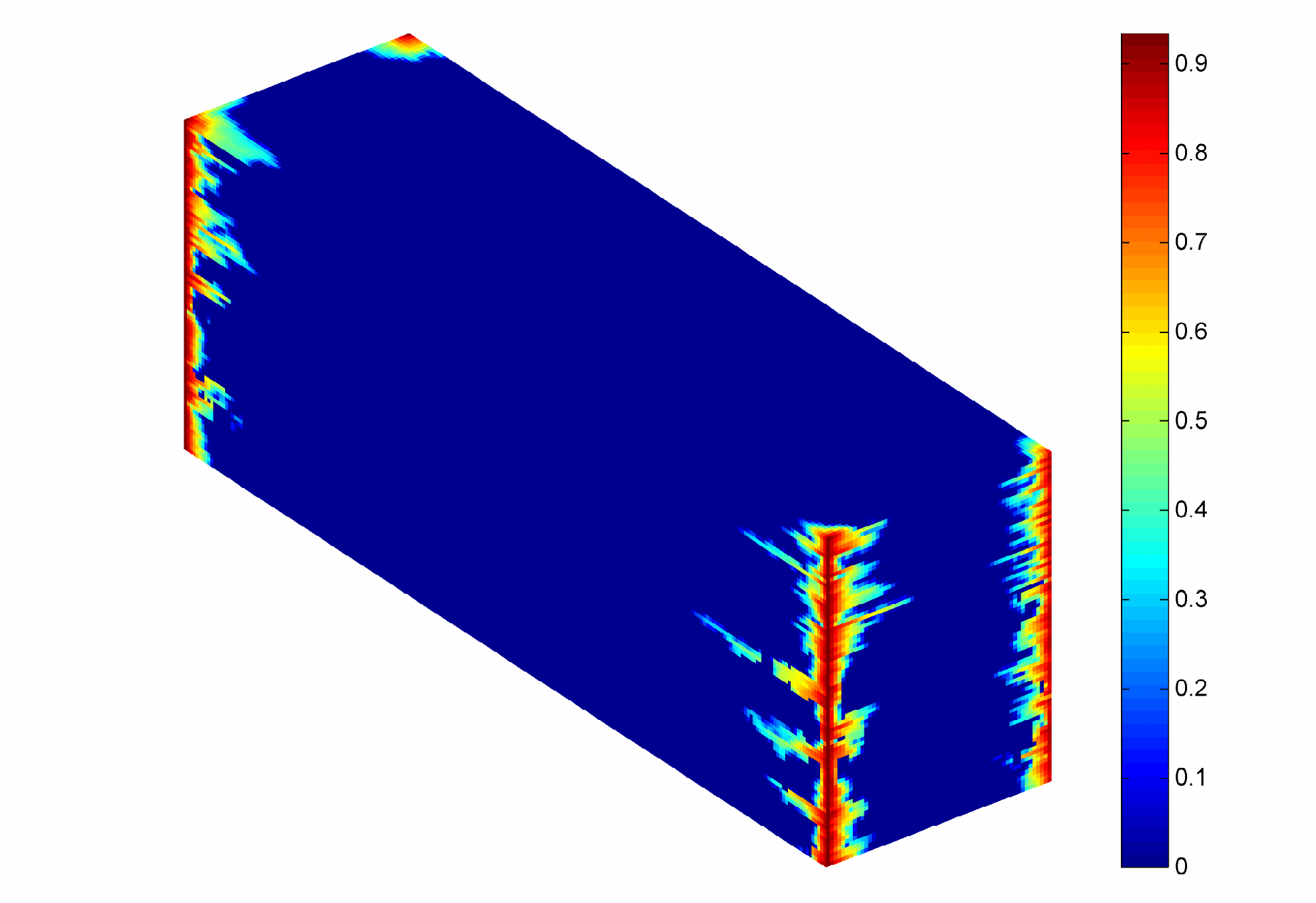}}
		\caption{Comparison of saturation obtained in three methods at time instant 50 for case 1.}
		\label{fig:scompare1}
	\end{figure}

\begin{figure}[H]
		\centering	
		\includegraphics[trim={1cm 6cm 1cm 7cm},clip,width=3.8in]{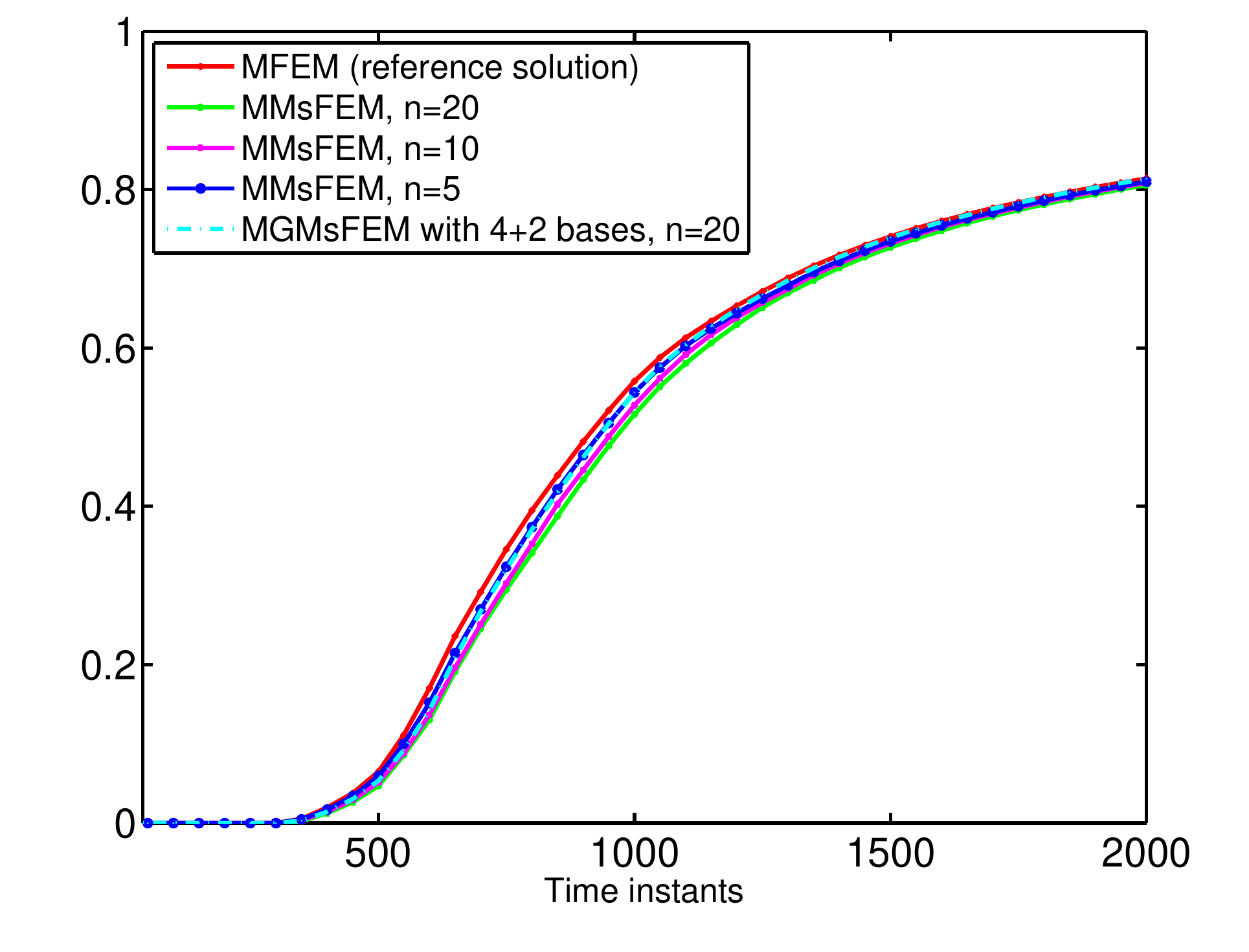}
		\caption{Comparison of water-cut between MMsFEM and MGMsFEM for case 2 at producer 1.}
		\label{fig:pc_b1}
	\end{figure}

\begin{figure}[H]
		\centering	
		\includegraphics[trim={1cm 6cm 1cm 7cm},clip,width=3.8in]{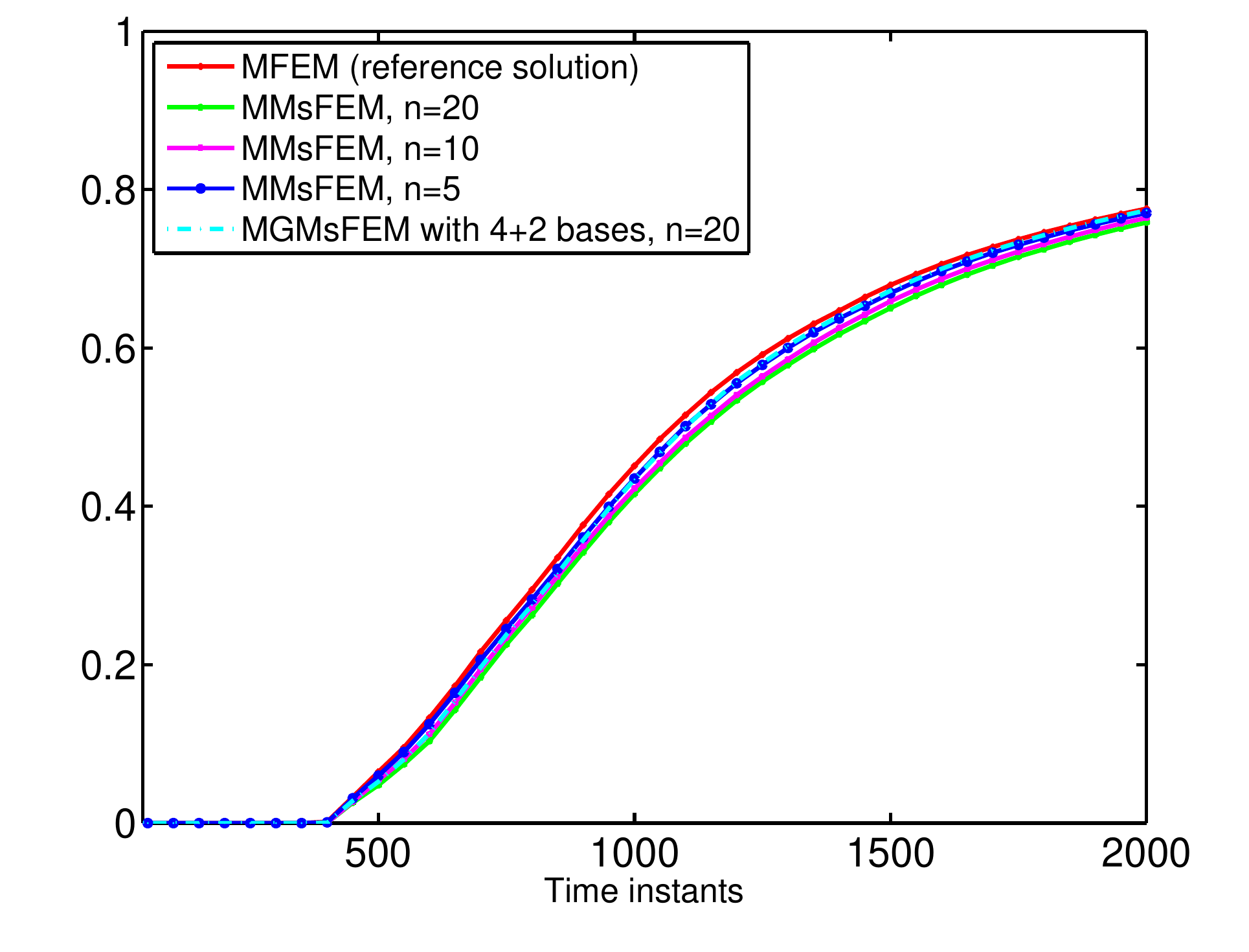}
		\caption{Comparison of water-cut between MMsFEM and MGMsFEM for case 2 at producer 2.}
		\label{fig:pc_b2}
	\end{figure}

\begin{figure}[H]
		\centering	
		\includegraphics[trim={1cm 6cm 1cm 7cm},clip,width=3.8in]{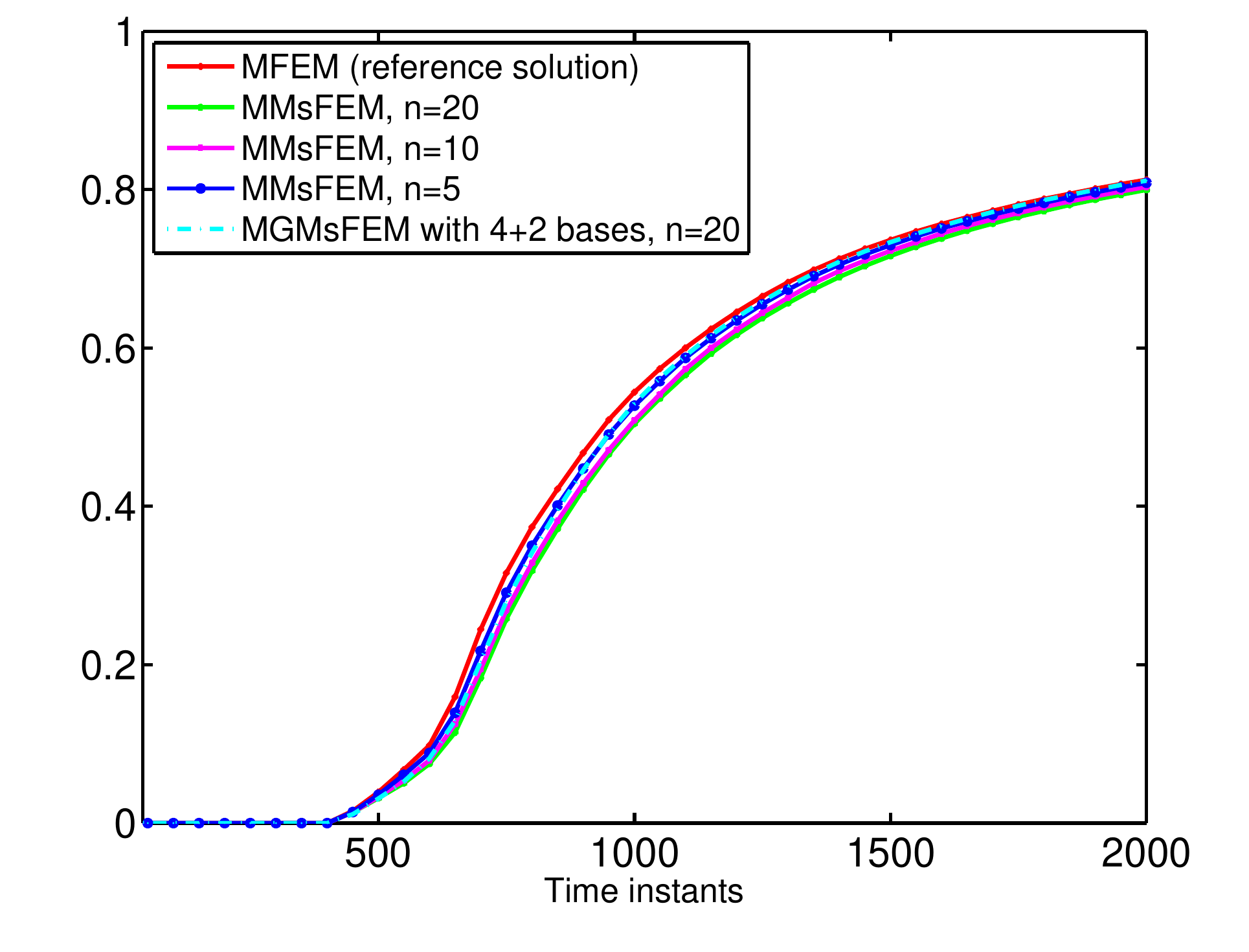}
		\caption{Comparison of water-cut between MMsFEM and MGMsFEM for case 2 at producer 3.}
		\label{fig:pc_b3}
	\end{figure}

\begin{figure}[H]
		\centering	
		\includegraphics[trim={1cm 6cm 1cm 7cm},clip,width=3.8in]{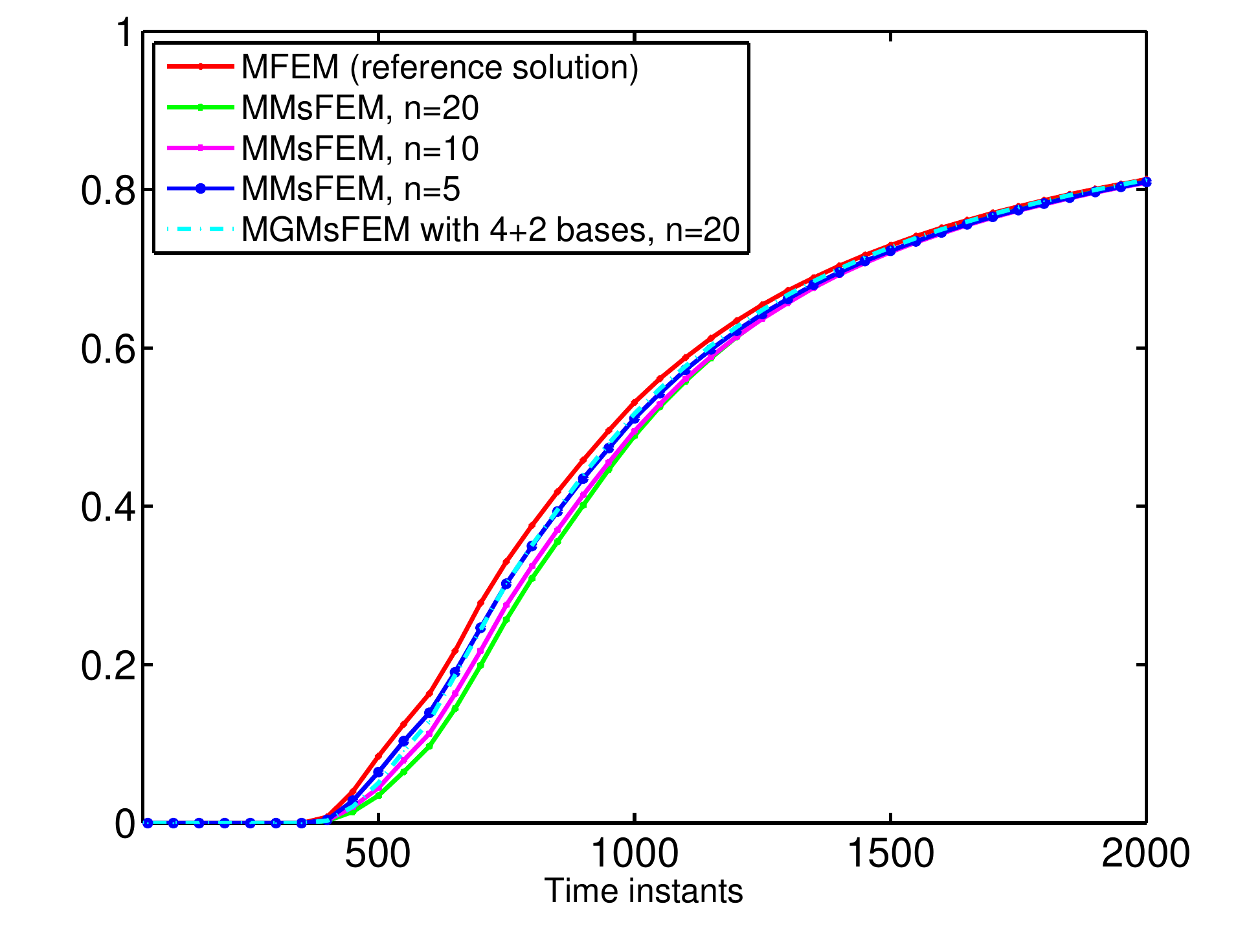}
		\caption{Comparison of water-cut between MMsFEM and MGMsFEM for case 2 at producer 4.}
		\label{fig:pc_b4}
	\end{figure}

	Figure \ref{fig:scompare1} and Figures \ref{fig:pc_b1}-\ref{fig:pc_b4} present the water-cut curves for case 1 and case 2 separately. There is only one producer while four in the latter case. For the water-cut, it is referred to the fraction of water $q_w/q_t$ in the produced fluid, where $q_t=q_w+q_o$. Besides, it is a function of time.  $q_w$ and $q_o$ are the flow rates of water and oil at the production edge of the model. Specifically, $q_w=\int_{\partial\Omega^{\text{out}}} f(S)v\cdot ndS$ and $q_t=\int_{\partial\Omega^{\text{out}}} v\cdot ndS$, where $\partial\Omega^{\text{out}}$ is the outer flow boundary. We make a comparison among four cases: MFEM, MMsFEM with n=5,10,20 and MGMsFEM with 4+2 basis and n=20. One may observe that the MGMsFEM best approximate the reference case while the MMsFEM offers relatively lower accuracy, where the accuracy increases with smaller coarse-mesh sizes. However, it is worth mentioning that the water-cut lines are almost identical between the best two approximations, MGMsFEM and MMsFEM with $n=5$.
	
\begin{figure}[H]
		\centering
		\subfigure[reference]{
			\includegraphics[trim={1cm 7.5cm 1cm 7cm},clip,width=3.0in]{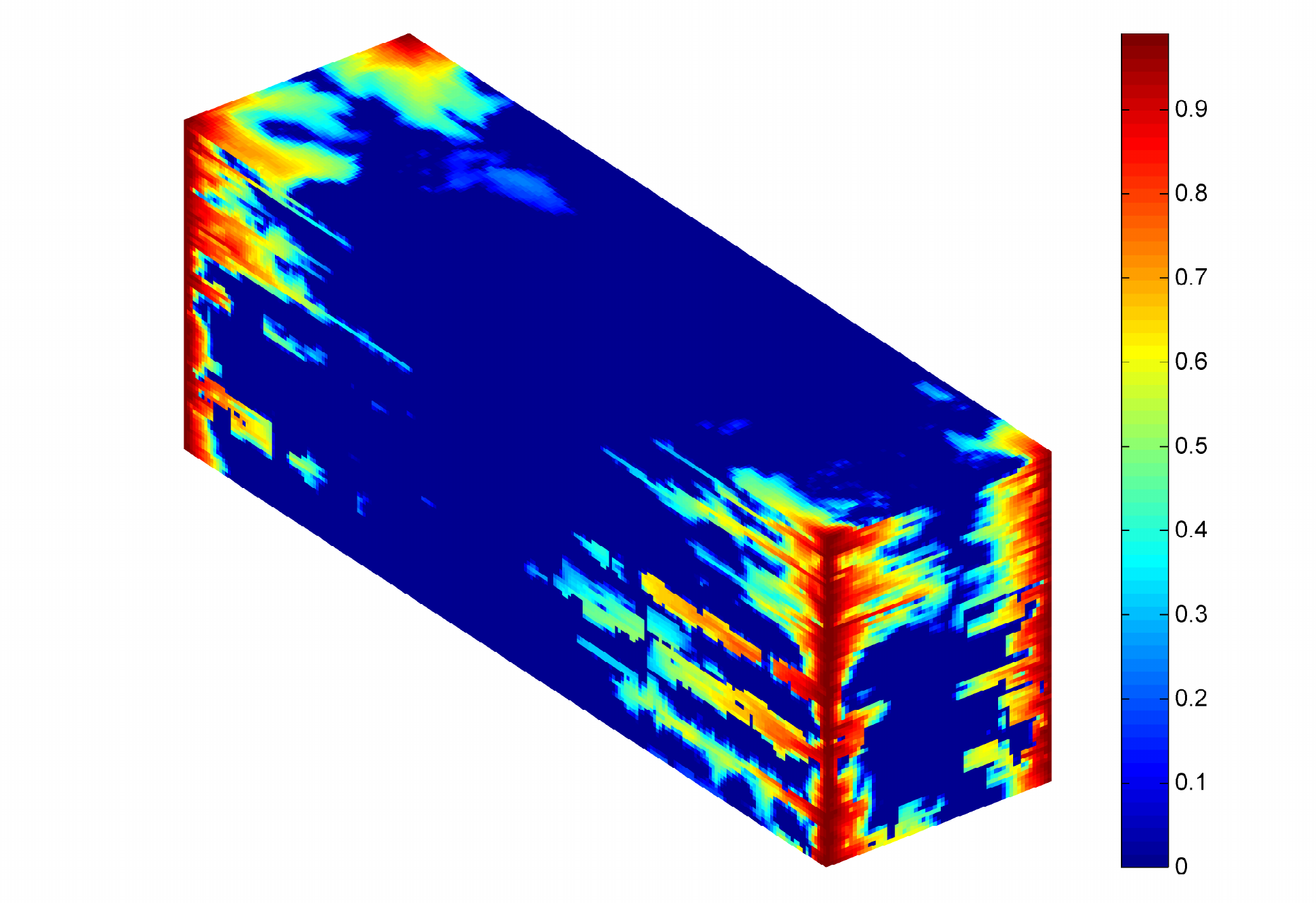}}
		\subfigure[MMsFEM with n=10]{
			\includegraphics[trim={1cm 7.5cm 1cm 7cm},clip,width=3.0in]{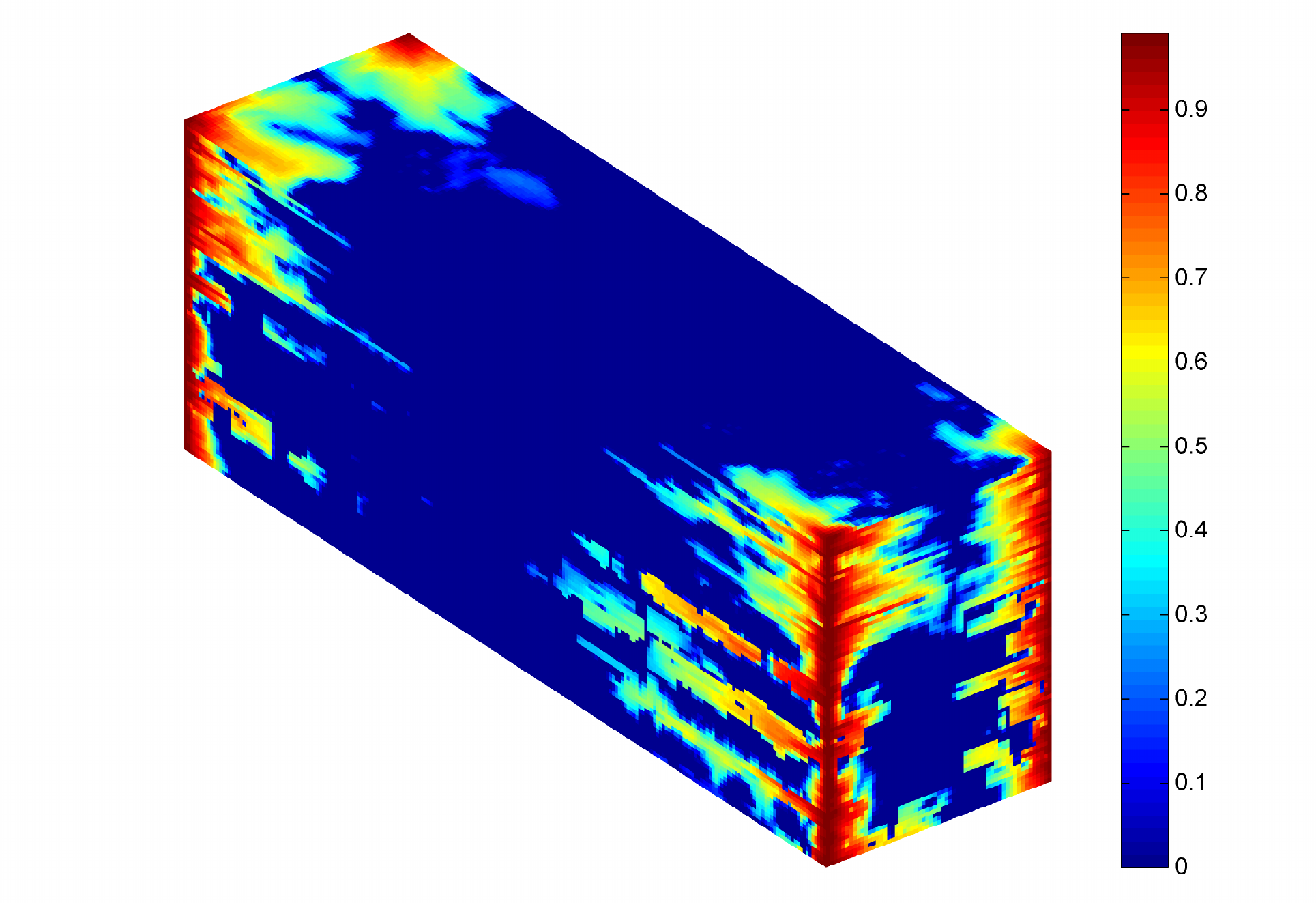}}
		\subfigure[MGMsFEM with n=20 and 4+2 basis functions]{
			\includegraphics[trim={1cm 7.5cm 1cm 7cm},clip,width=3.0in]{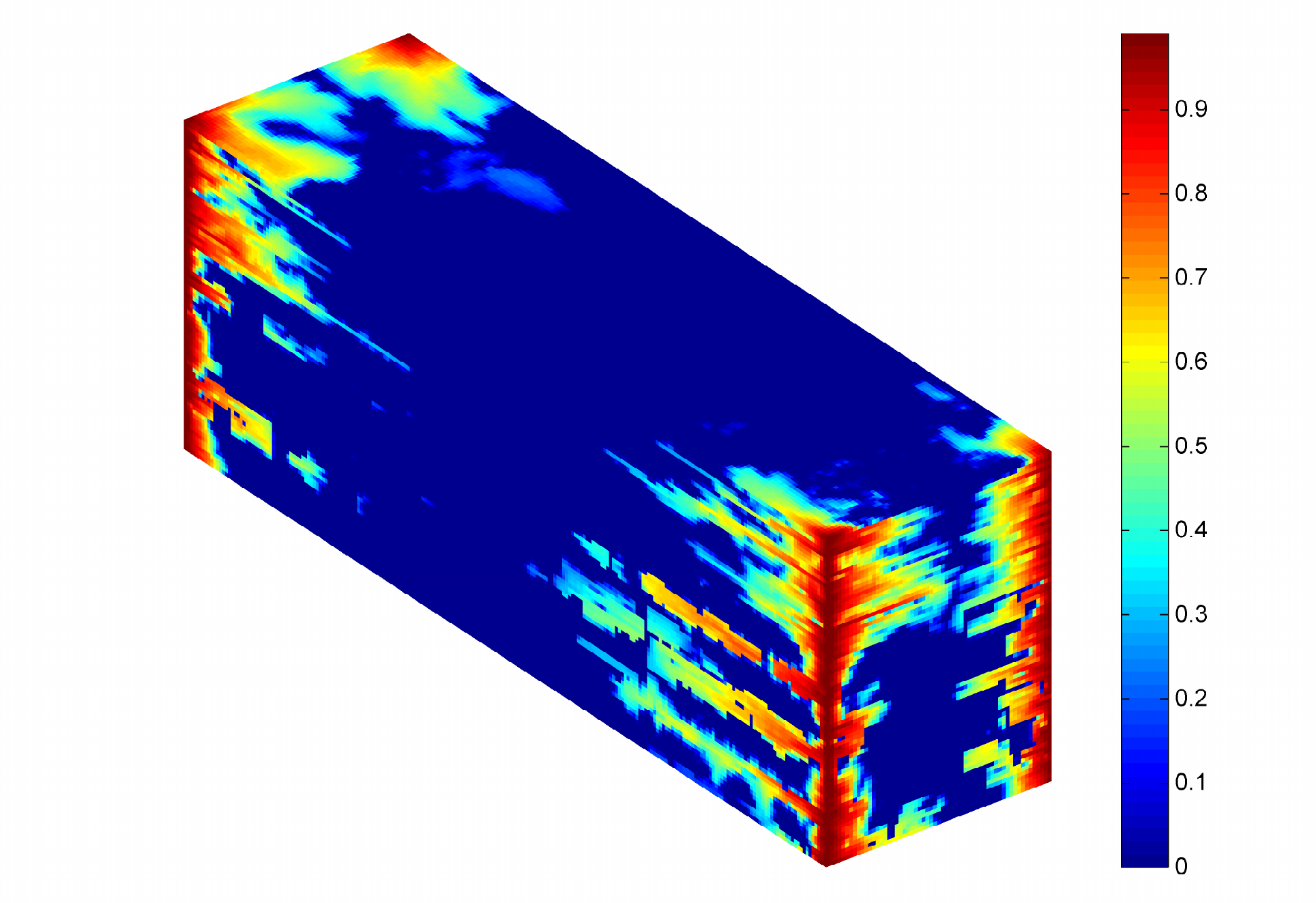}}
		\caption{Comparison of saturation obtained in three methods at time instant 500 for case 1.}
		\label{fig:scompare2}
	\end{figure}

\begin{figure}[H]
		\centering
		\subfigure[reference]{
			\includegraphics[trim={1cm 7.5cm 1cm 7cm},clip,width=3.0in]{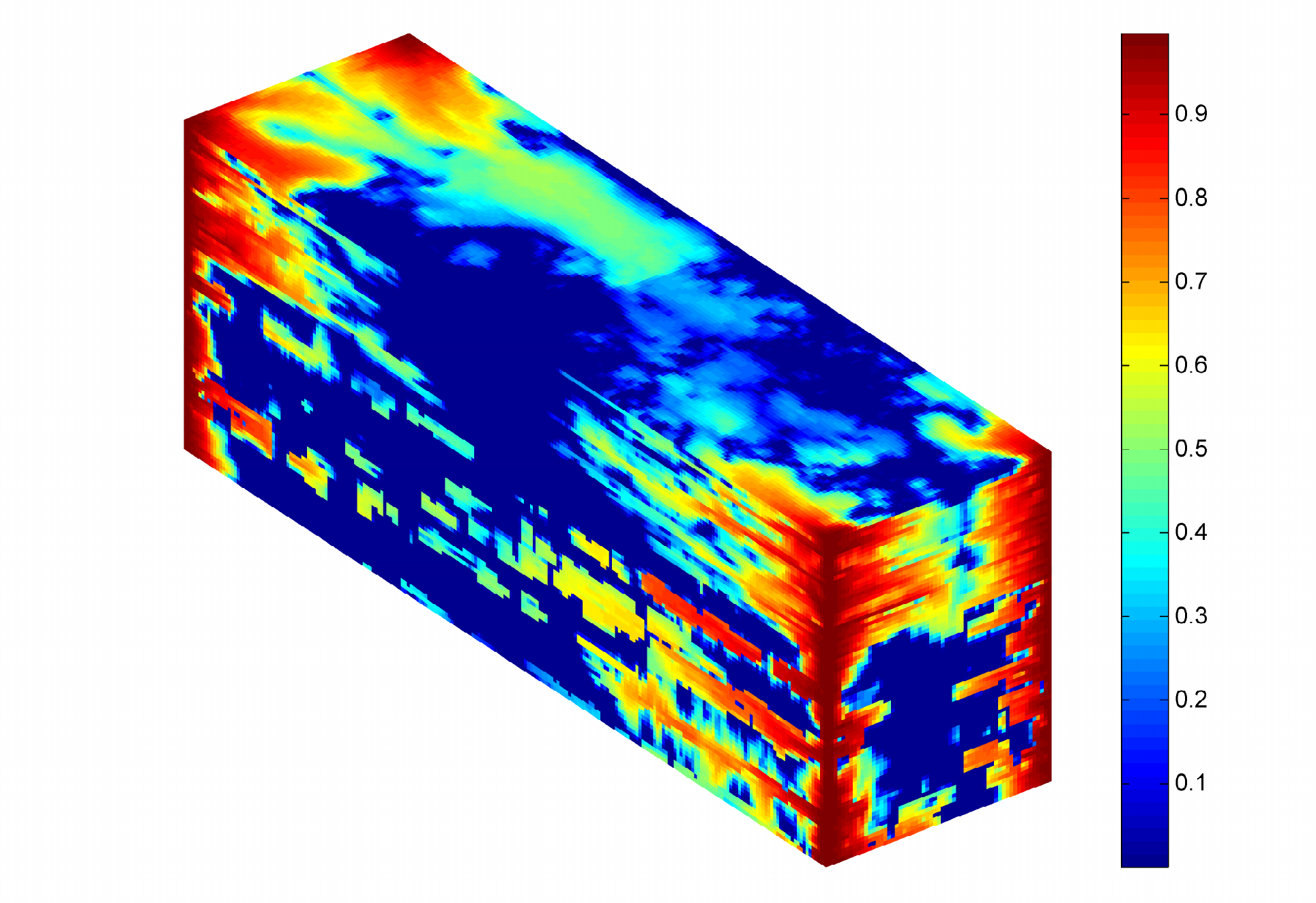}}
		\subfigure[MMsFEM with n=10]{
			\includegraphics[trim={1cm 7.5cm 1cm 7cm},clip,width=3.0in]{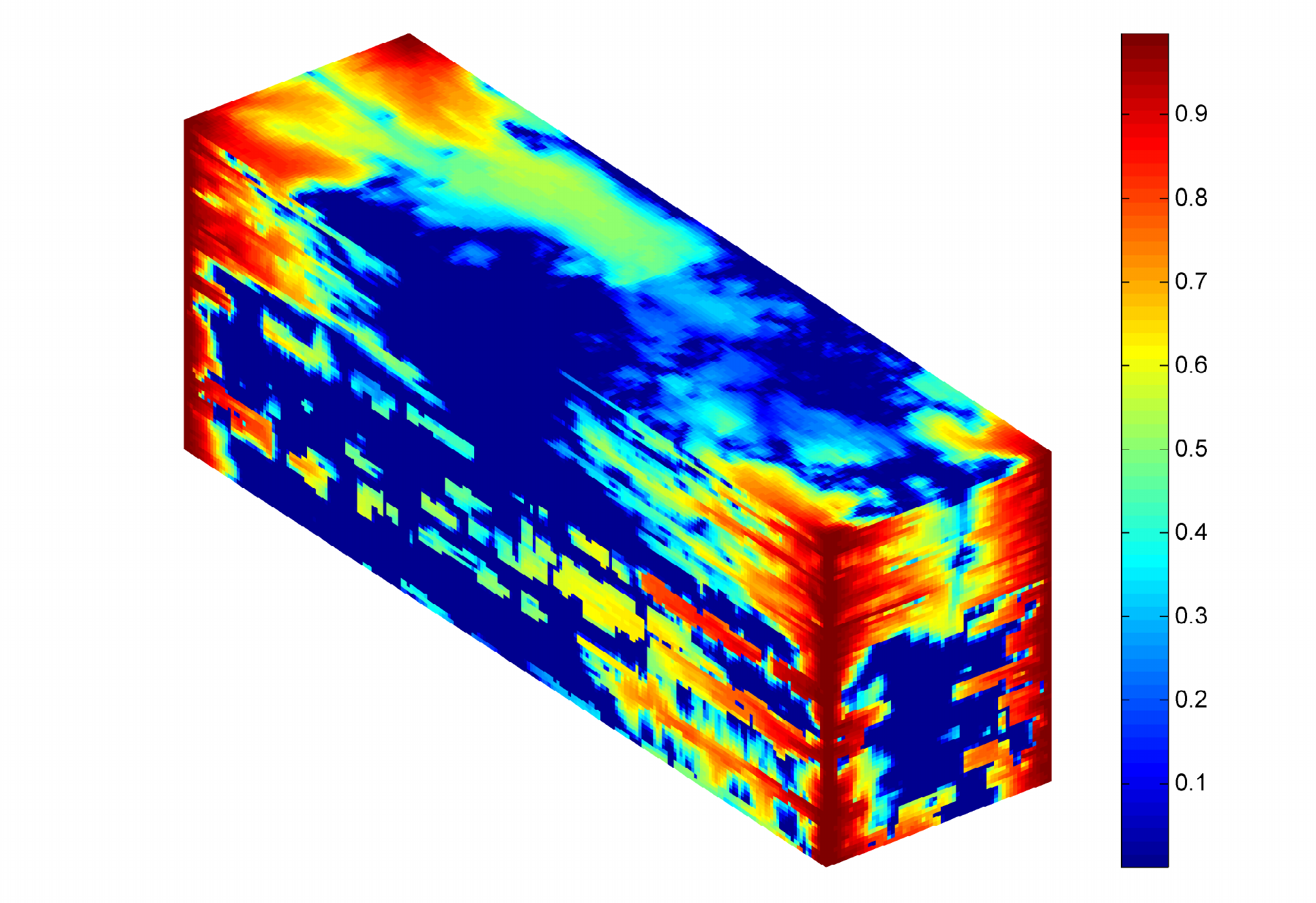}}
		\subfigure[MGMsFEM with n=20 and 4+2 basis functions]{
			\includegraphics[trim={1cm 7.5cm 1cm 7cm},clip,width=3.0in]{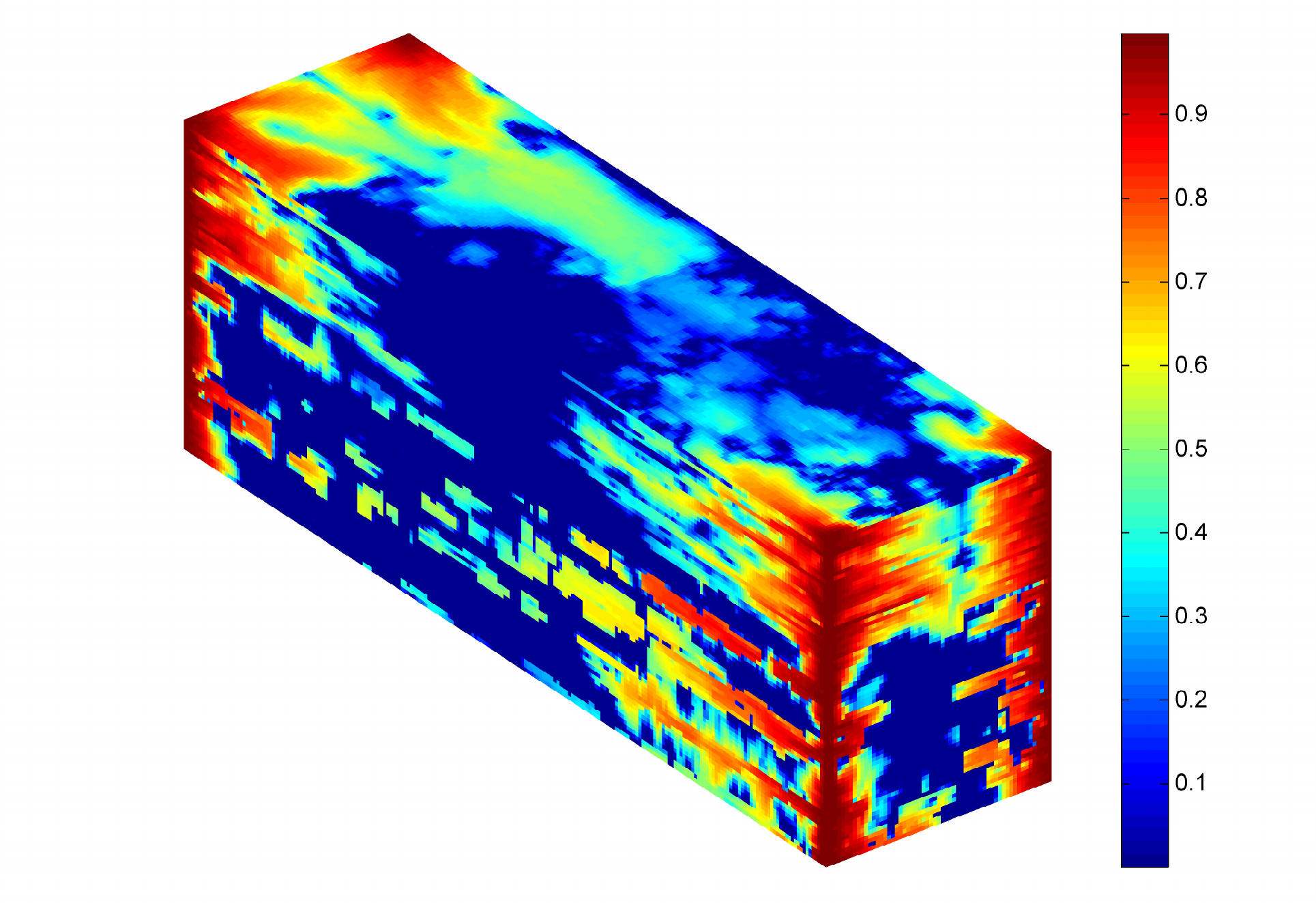}}
		\caption{Comparison of saturation obtained in three methods at time instant 2000 for case 1}
		\label{fig:scompare3}
	\end{figure}

\begin{figure}[H]
		\centering
		\subfigure[reference]{
			\includegraphics[trim={1cm 7.5cm 1cm 7cm},clip,width=3.0in]{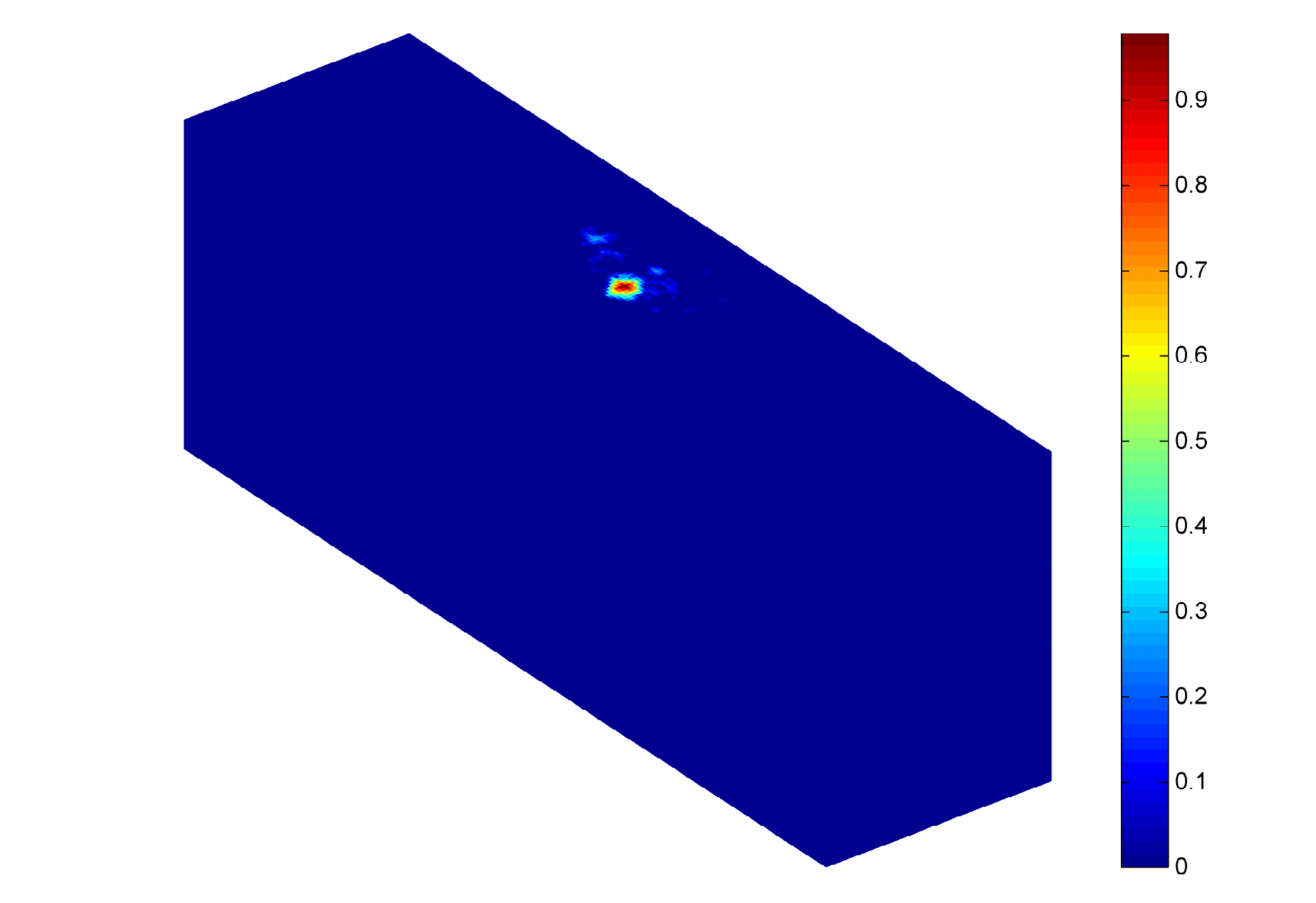}}
		\subfigure[MMsFEM with n=10]{
			\includegraphics[trim={1cm 7.5cm 1cm 7cm},clip,width=3.0in]{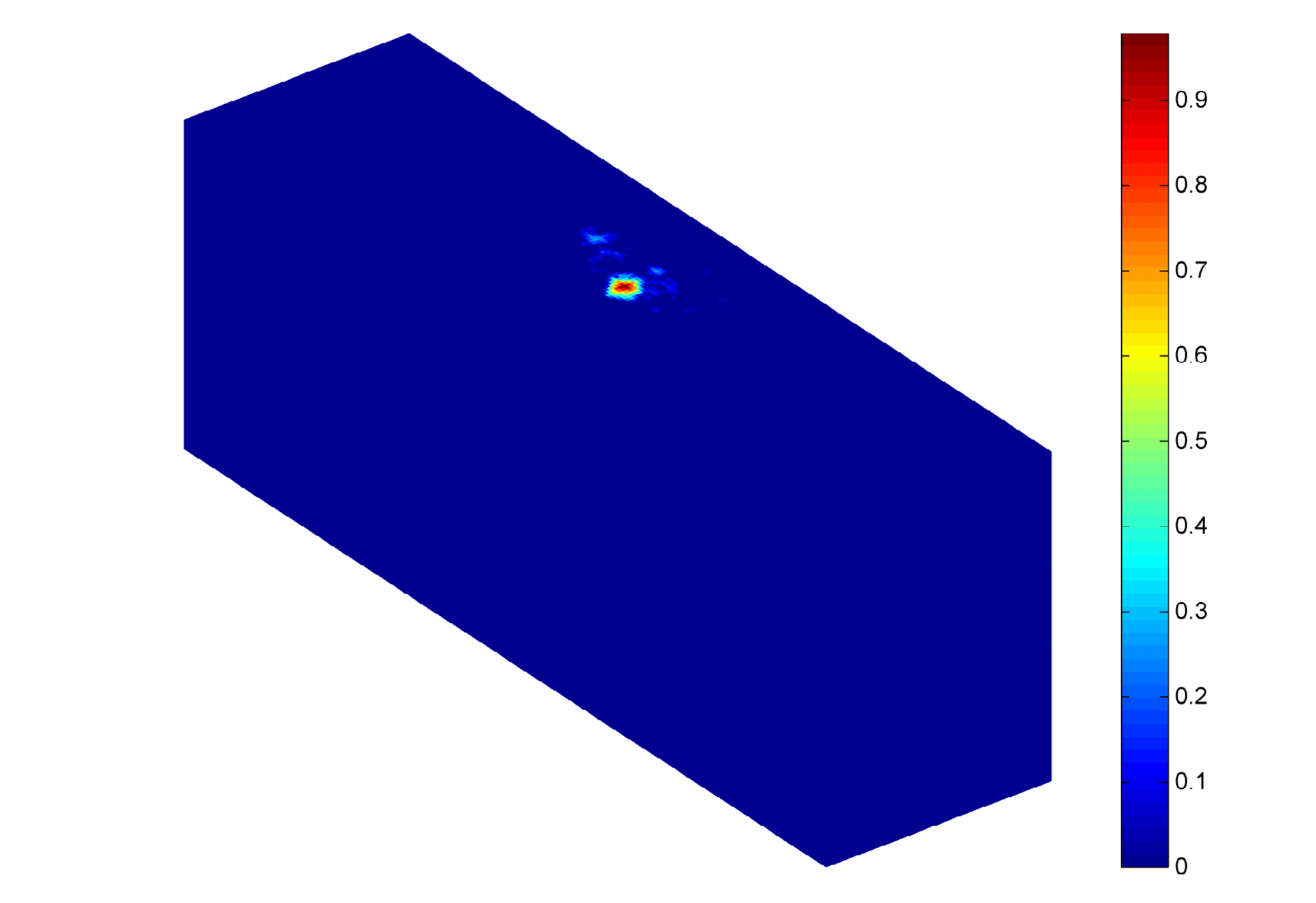}}
		\subfigure[MGMsFEM with n=20 and 4+2 basis functions]{
			\includegraphics[trim={1cm 7.5cm 1cm 7cm},clip,width=3.0in]{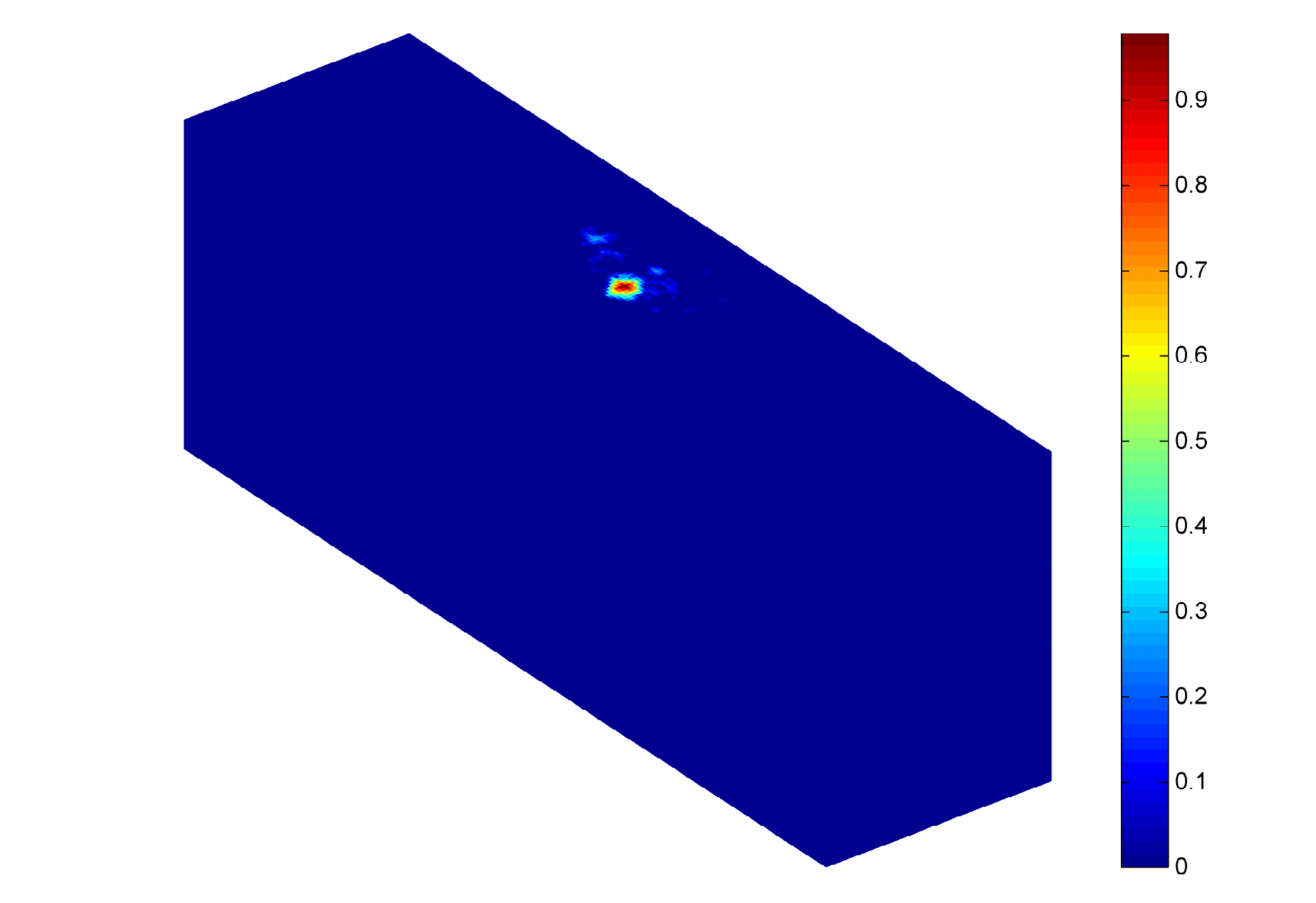}}
		\caption{Comparison of saturation obtained in three methods at time instant 50 for case 2.}
		\label{fig:scompareb1}
	\end{figure}

\begin{figure}[H]
		\centering
		\subfigure[reference]{
			\includegraphics[trim={1cm 7.5cm 1cm 7cm},clip,width=3.0in]{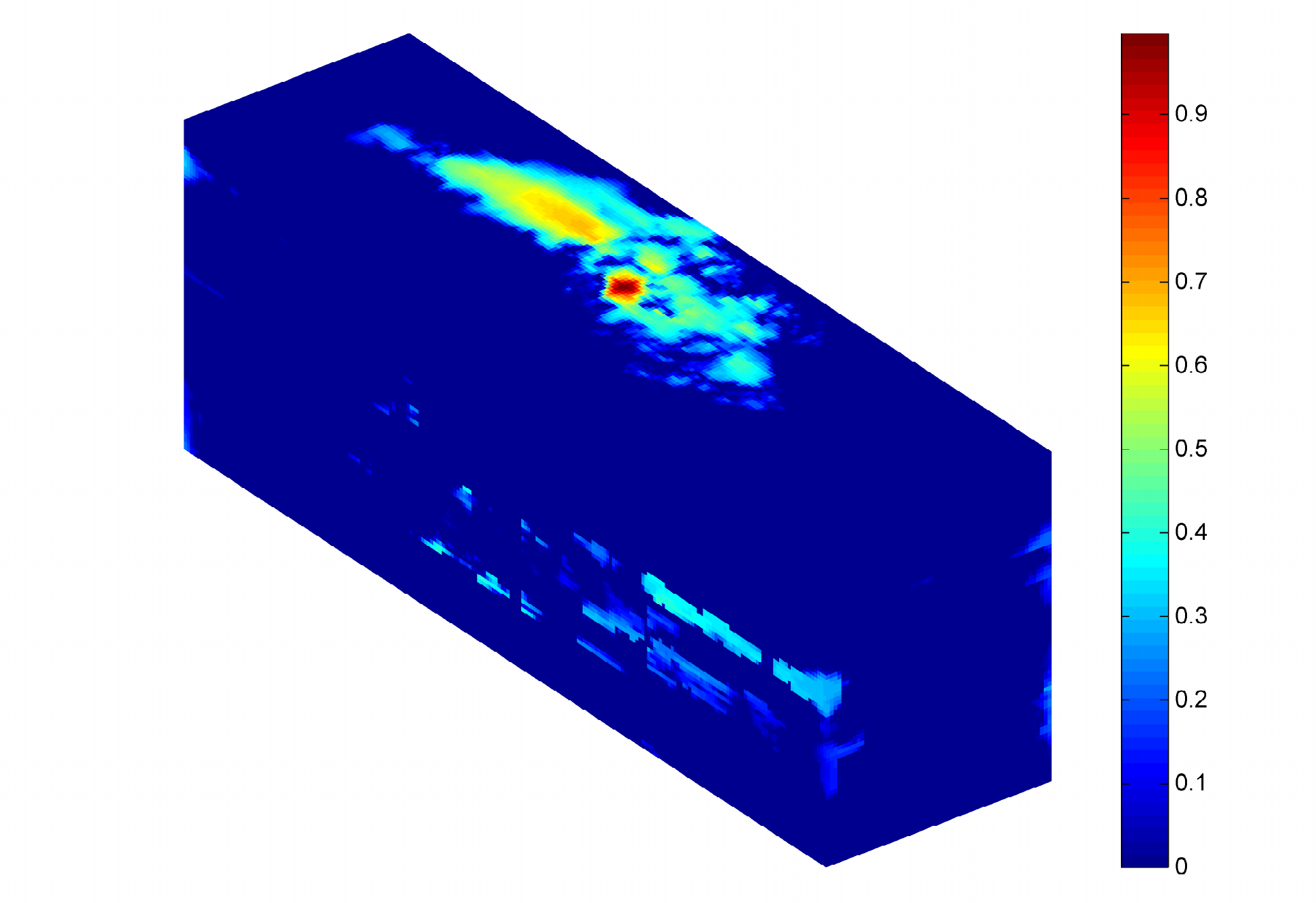}}
		\subfigure[MMsFEM with n=10]{
			\includegraphics[trim={1cm 7.5cm 1cm 7cm},clip,width=3.0in]{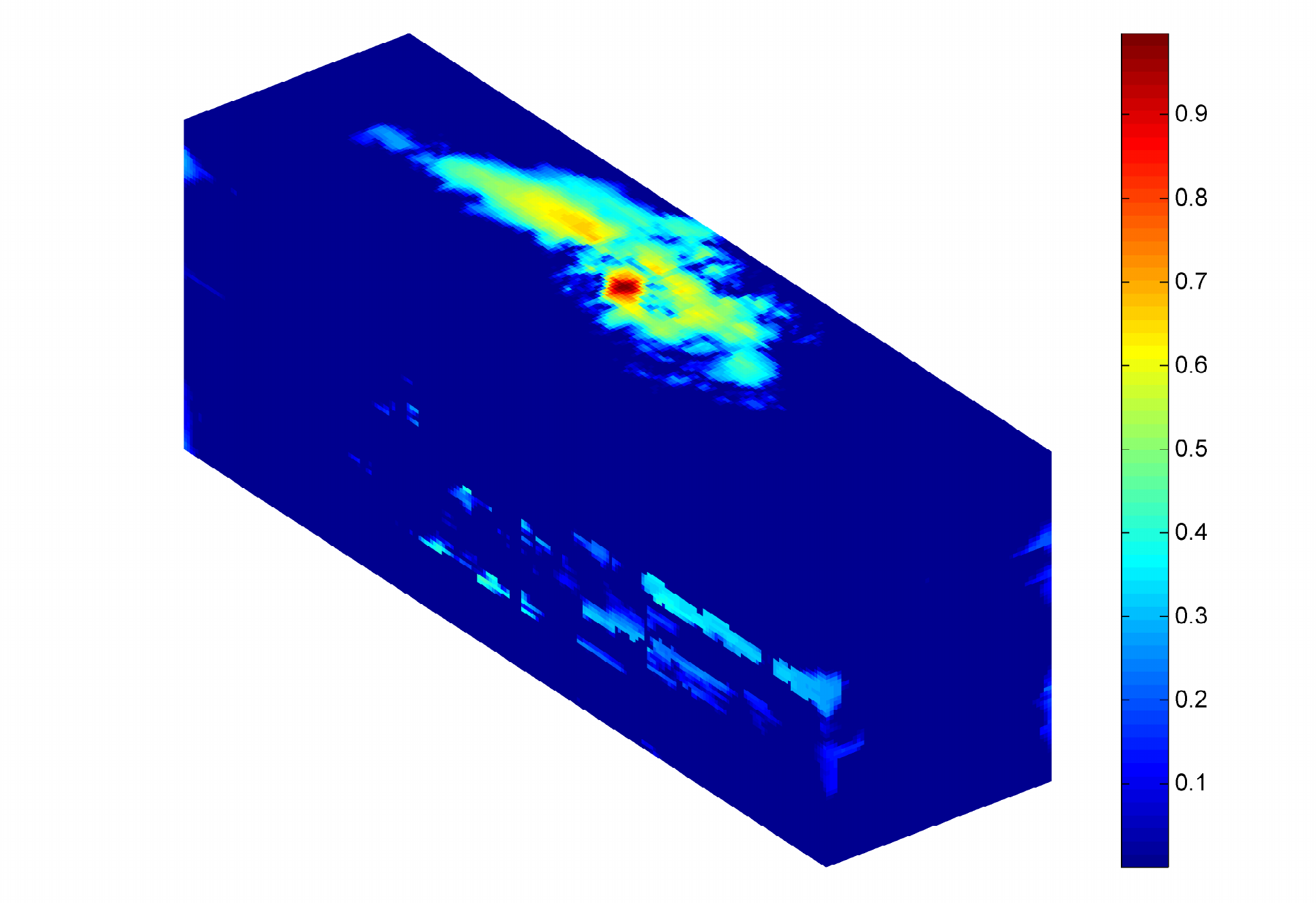}}
		\subfigure[MGMsFEM with n=20 and 4+2 basis functions]{
			\includegraphics[trim={1cm 7.5cm 1cm 7cm},clip,width=3.0in]{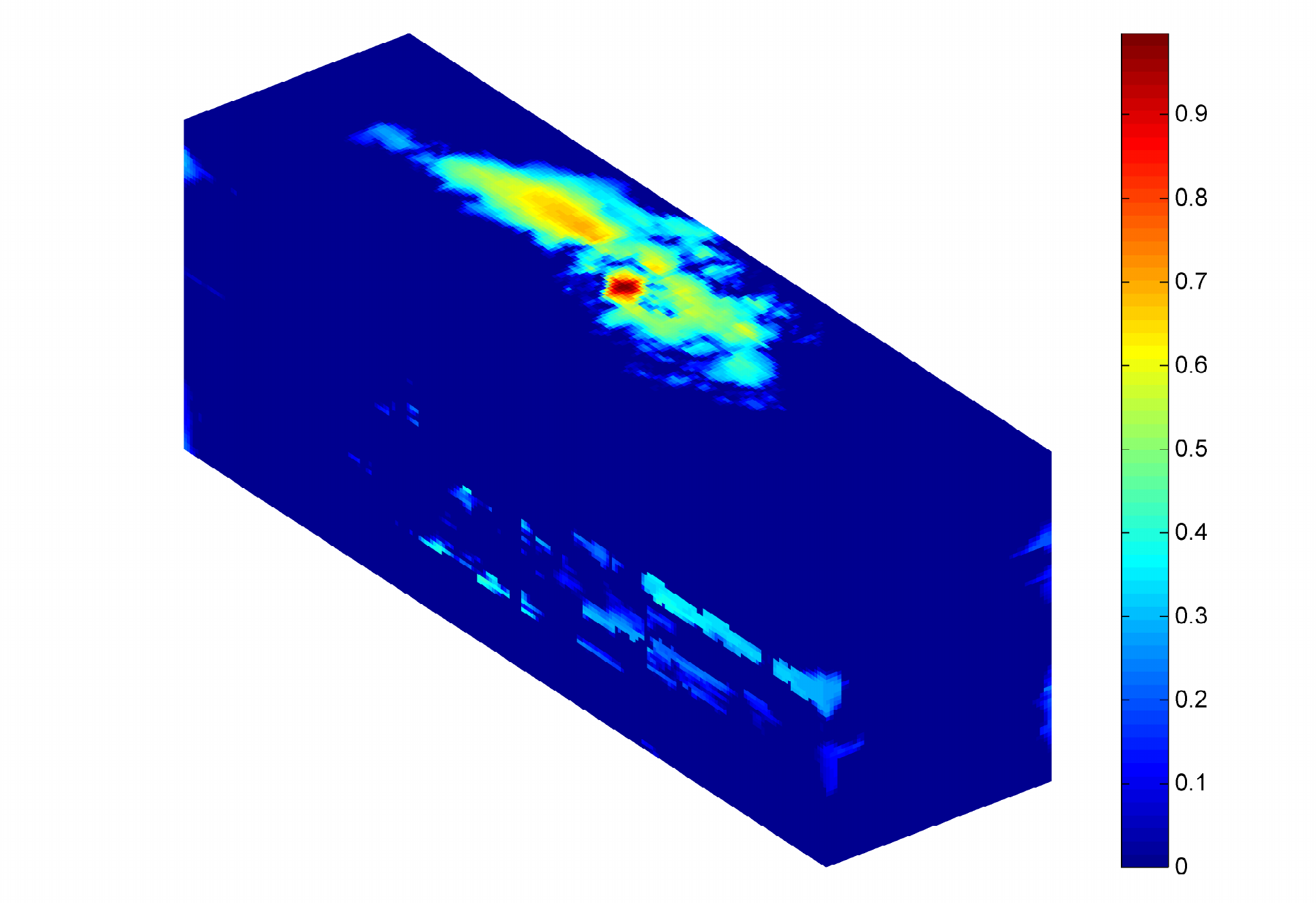}}
		\caption{Comparison of saturation obtained in three methods at time instant 500 for case 2.}
		\label{fig:scompareb2}
	\end{figure}

\begin{figure}[H]
		\centering
		\subfigure[reference]{
			\includegraphics[trim={1cm 7.5cm 1cm 7cm},clip,width=3.0in]{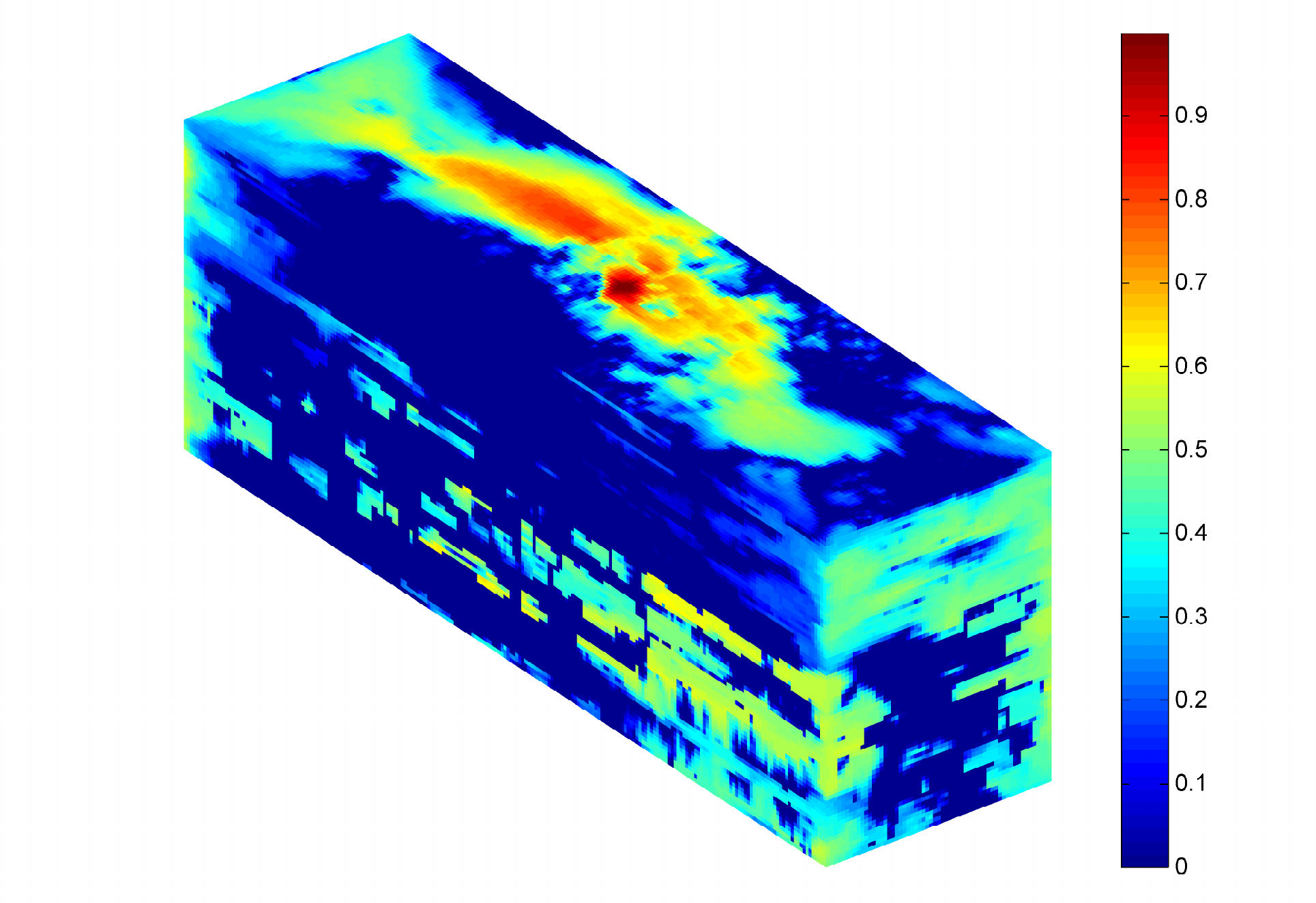}}
		\subfigure[MMsFEM with n=10]{
			\includegraphics[trim={1cm 7.5cm 1cm 7cm},clip,width=3.0in]{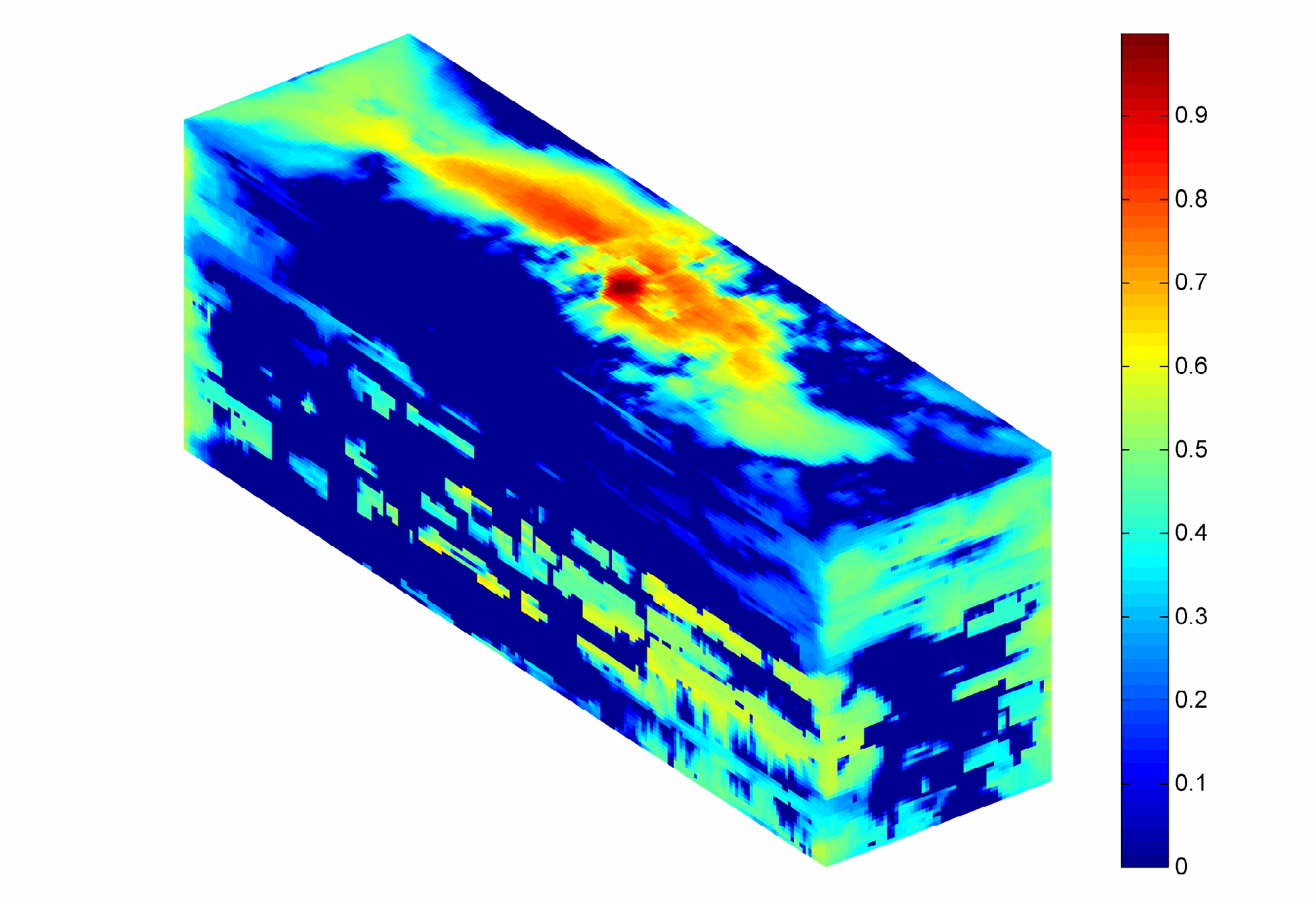}}
		\subfigure[MGMsFEM with n=20 and 4+2 basis functions]{
			\includegraphics[trim={1cm 7.5cm 1cm 7cm},clip,width=3.0in]{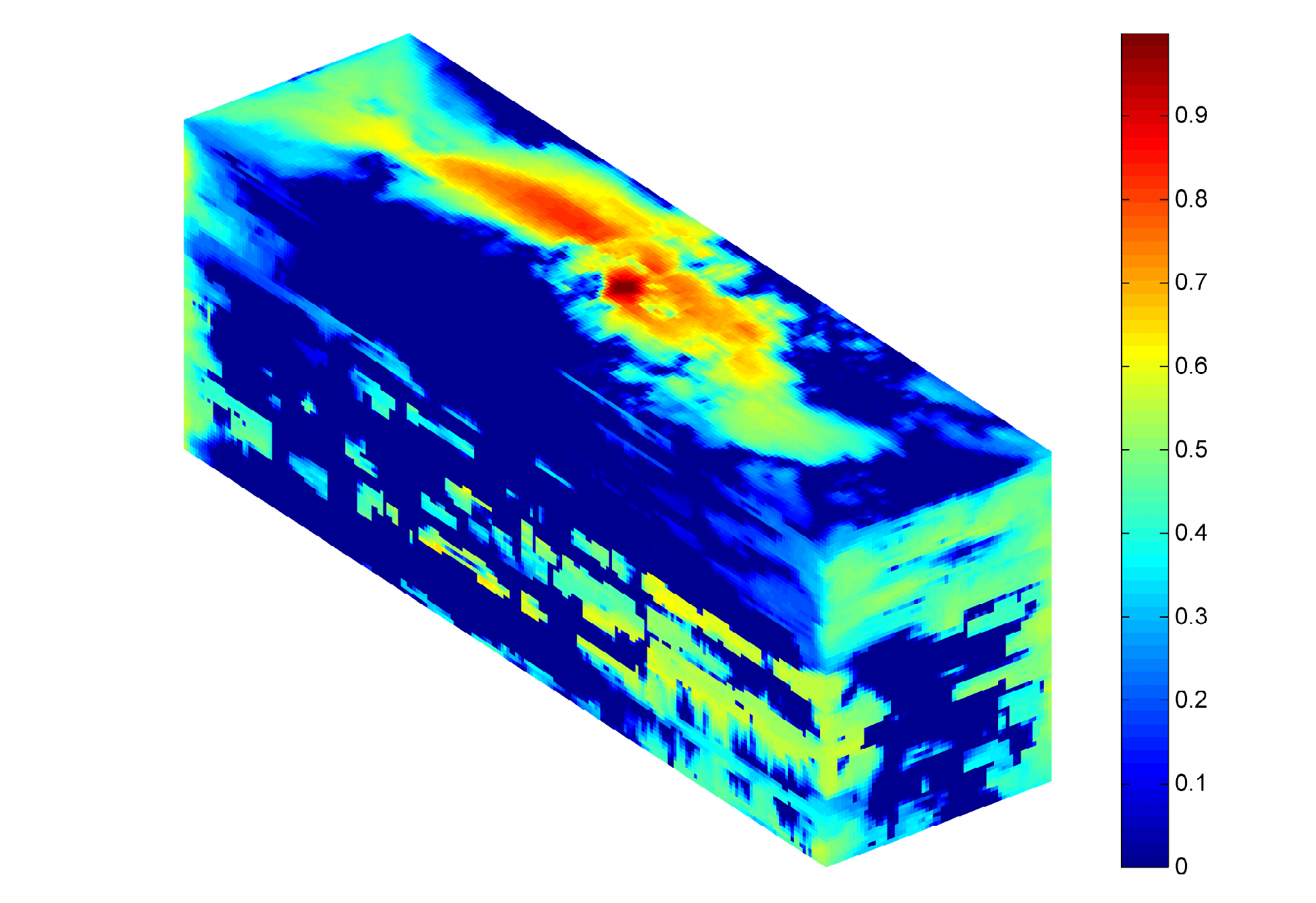}}
		\caption{Comparison of saturation obtained in three methods at time instant 2000 for case 2.}
		\label{fig:scompareb3}
	\end{figure}

	Figure \ref{fig:scompare1}-\ref{fig:scompare3} show saturation dynamics at three time instants, 50, 500, 2000 respectively, where water is injected from four sides. In each figure, (a),(b),(c) correspond to solutions resulted from MFEM, MMsFEM with $n=10$ and MGMsFEM with $n=20$ and $4+2$ basis in each local neighborhood. It is evident that solution by MGMsFEM can approximate the reference case excellently while some minor differences could be detected from MFEM and MMsFEM. However, in addition to the fact that the overall shapes of saturation maps are similar between reference and two approximations, most of the fine detail is retained, which shows the computational efficiency of these two approximation methods. Figure \ref{fig:scompareb1}-\ref{fig:scompareb3} show the dynamics of saturation in the other circumstance, where the injector is in the center. Similarly results can be attained in this case. Consequently, these two methods offer relatively good alternatives for reference solution from MFEM, which is consistent with aforementioned error results.  Summarizing these cases, numerical results convincingly show that one can use single-phase flow computation to
	solve the two-phase flow model to get relatively accurate approximation.\\
	
	\begin{table}[H]
		\centering \begin{tabular}{|c|c|c|c|c|c|c|c|c|}
			\hline
			Method    & $n$ &Dof   & $T_{\text{setup}}(s)$ &$T_{\text{sim}}(s)$ &$e_s$  \tabularnewline
			\hline
			MFEM&/   &4188400&/&122879.0 &/ \tabularnewline\hline
			MMsFEM&20    &439 &3114.4& 3363.4&0.091  \tabularnewline\hline
			MMsFEM&10   &3868 &3338.7&6570.9 &0.084  \tabularnewline\hline
			MMsFEM&5    &32368&4128.0&28505.6 &0.060  \tabularnewline\hline
			MGMsFEM(1+0)&20    & 439&1514.5&3823.4 &0.393 \tabularnewline\hline
			MGMsFEM(1+0)&10    & 3868&625.9&7966.3 &0.311 \tabularnewline\hline
			MGMsFEM(1+0)&5    & 32368&3093.4&31715.5 &0.168 \tabularnewline\hline
			MGMsFEM(6+0)&20    &1974&1669.6&5451.3 &0.103 \tabularnewline\hline
			MGMsFEM(8+0)&20   & 2588&1718.4&6059.4 &0.080 \tabularnewline\hline
			MGMsFEM(2+1)&20    &1053 &2846.4&3741.4 & 0.112\tabularnewline\hline
			MGMsFEM(2+2)&20    & 1360&3992.2&4175.7 & 0.064\tabularnewline\hline
			MGMsFEM(4+2)&20   & 1974&4118.8&4997.4 &0.047 \tabularnewline\hline
			MGMsFEM(2+2)&10    & 12304&1767.0&8258.0 &0.046 \tabularnewline\hline
			
		\end{tabular}
		\caption{Computation time and error comparison of different methods for case 1.}
		\label{time1}
	\end{table}

	\begin{table}[H]
		\centering \begin{tabular}{|c|c|c|c|c|c|c|c|c|}
			\hline
			Method    & $n$ &Dof& $T_{\text{setup}}(s)$ &$T_{\text{sim}}(s)$ &$e_s$  \tabularnewline
			\hline
			MFEM&/ &4188400 &/&119944.2 &/ \tabularnewline\hline		
			MMsFEM&20  &439 &2930.3&3756.7 &0.096  \tabularnewline\hline
			MMsFEM&10  &3868& 3069.3&6689.9 &0.082  \tabularnewline\hline
			MMsFEM&5   &32368& 4180.2&31195.6 &0.063  \tabularnewline\hline
			MGMsFEM(1+0)&20   &439&1539.1&4126.2&0.353  \tabularnewline\hline
			MGMsFEM(1+0)&10   &3868& 637.6&7971.9&0.297  \tabularnewline\hline
			MGMsFEM(1+0)&5   &32368& 3146.1&32127.2&0.167  \tabularnewline\hline
			MGMsFEM(6+0)&20    &1974&1708.9&5108.2 &0.110 \tabularnewline\hline
			MGMsFEM(8+0)&20    &2588&1727.8&5961.6 &0.086 \tabularnewline\hline	
			MGMsFEM(2+1)&20    &1053&2884.2&4017.7 &0.122 \tabularnewline\hline		
			MGMsFEM(2+2)&20    &1360&4076.4&4540.6 &0.077 \tabularnewline\hline
			MGMsFEM(4+2)&20    &1974&4178.8&5203.4 & 0.059\tabularnewline\hline
			MGMsFEM(2+2)&10    &12304&1771.3&8851.4 &0.051 \tabularnewline\hline		
		\end{tabular}
		\caption{Computation time and error comparison of different methods for case 2.}
		\label{time2}
	\end{table}
	\begin{figure}[H]
		\centering	
		\includegraphics[trim={1cm 6cm 1cm 7cm},clip,width=3.8in]{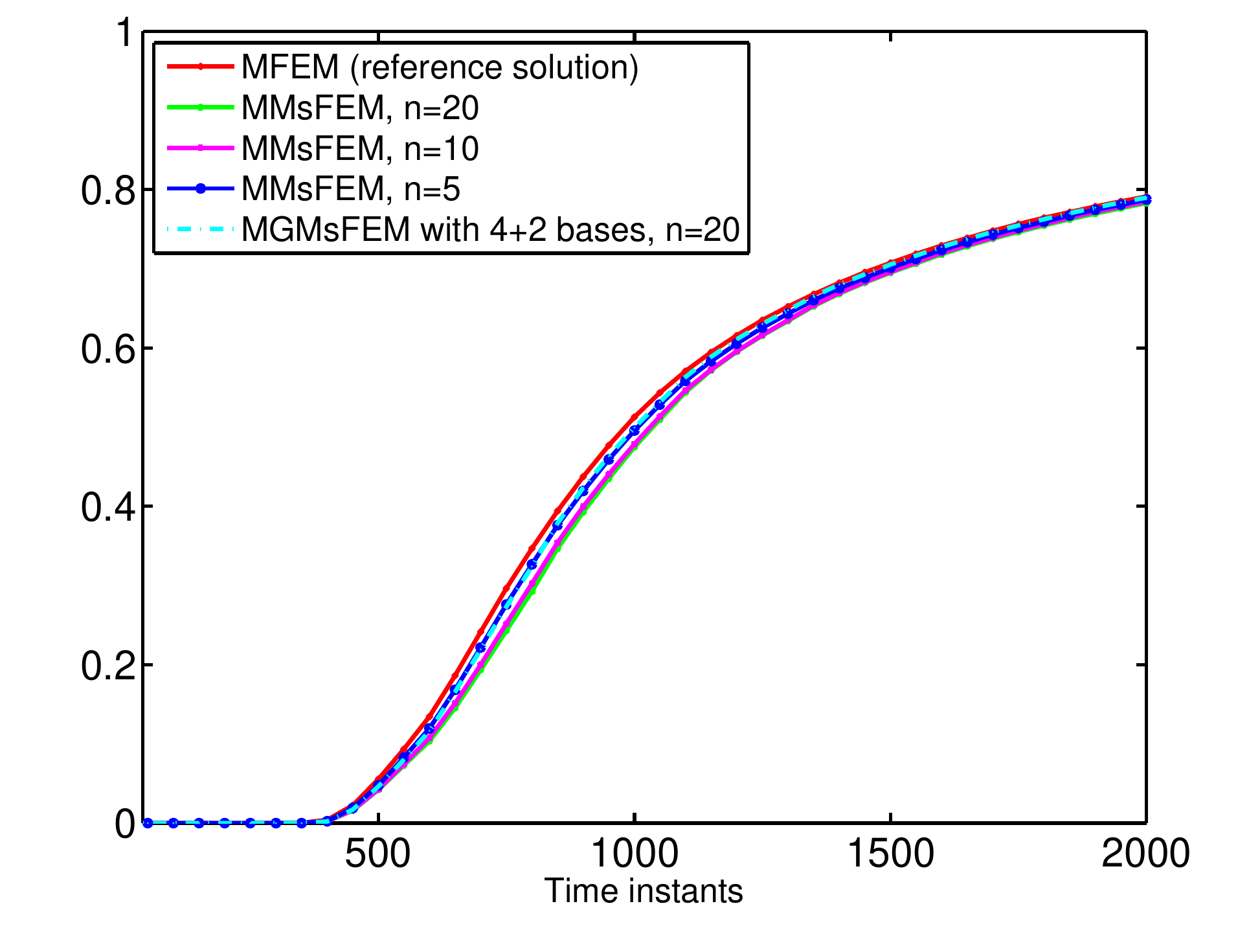}
		\caption{Comparison of water-cut between MMsFEM and MGMsFEM for case 1.}
		\label{fig:pc_1}
	\end{figure}
	
	
	\section{Discussions and conclusions}
	In the paper, we compared two multiscale methods for the simulations of two-phase models with permeability field which is highly heterogeneous and of high contrast. The main difference between these two methods is the construction of multiscale basis functions. For MGMsFEM, we construct the offline multiscale space via a set of independent spectral problems, which can be implemented in parallel to save time. For the online multiscale space, we enrich the offline space by some residual driven online basis functions, which may be solved in some oversampled regions. The global information contained in the online basis functions enhances the accuracy remarkably. For the MMsFEM with limited global information, it is modified from the standard MMsFEM. In this method, no local spectral problem is implemented. In addition, for MMsFEM, we incorporate global information in the computation of multiscale basis functions. In particular, we need to solve an independent single-phase problem in advance to get information about the global flow pattern from the initial velocity field. This procedure is not included in MGMsFEM but it serves as an effective tool to improve the performance of the MMsFEM.
	
	In terms of the most important factor, i.e. accuracy, we can see that these two methods can both give satisfactory outcomes. MMsFEM can retain most of the small-scale information from the reference solutions. Besides, it can give relatively good approximations of the concerned quantities, such as the water-cut and saturation. Another advantage of this method is that the multiscale space is smaller that that of MGMsFEM since only one multiscale basis function is incorporated in each local neighborhood, which results in higher efficiency. However, a shortcoming is that the performance of the method is dependent on the similarity of the initial velocity field and the actual one in simulation. If there is a big mismatch between these two cases, there is no guarantee that the MMsFEM will give a satisfying result. Moreover, since single basis function is required in each local neighborhood for MMsFEM, we need to seek other ways to improve the performance. On the other hand, there are more flexibility for the MGMsFEM since we can change the number of offline and online basis functions in each local region. Besides, we can also adjust the number of oversampling layers as we construct the online basis functions. Consequently, MGMsFEM is a more flexible method than the MMsFEM and one may easily find a combination of the concerned parameters to obtain a better approximation. However, there are still some drawbacks for MGMsFEM. First, as one may verify from the numerical results, in somes cases only offline multiscale space is not sufficient to obtain a good approximation especially when the permeability field is of high heterogeneity such as $\kappa_3$ in this paper. In our experiments, MMsFEM can offer better results than MGMsFEM with a much bigger offline space. Consequently, one may seek to enrich the offline space with online basis functions, which require further computation. Besides, since multiple basis functions are used in each local neighborhood, the final system is bigger than MMsFEM, which give rise to larger computational cost. One also needs to construct a snapshot space and reduce its dimension to obtain the offline space, which is not required in MMsFEM.
	
	Consequently, the choice between these two methods is dependent on the actual problem. For example, if it is easy to obtain the required information of the initial velocity field or it is well-prepared, then MMsFEM may be considered as a good alternative. On the other hand, if it is demanding to obtain the required global information and offline space in MGMsFEM is sufficient to give a satisfying result, then MGMsFEM may be advantageous over the other one. However, when MMsFEM can not provide sufficient accuracy, one may apply MGMsFEM combined with online enrichment if the offline space is not sufficient. Thus, we can see that both MMsFEM and MGMsFEM possess some advantages, which lead to remarkable performances in some specific cases. However, both methods have some limitations. It is worth mentioning that in the following work we can combine these two methods in some way to pursue better outcomes than applying a single method.

	\section*{Acknowledgments}
	
	The research of Eric Chung is partially supported by the Hong Kong RGC General Research Fund (Project numbers 14304719 and 14302018) and CUHK Faculty of Science Direct Grant 2019-20.


\begin{thebibliography}{10}

\bibitem{aarnes2004use}
Jorg~E Aarnes.
\newblock On the use of a mixed multiscale finite element method for
  greaterflexibility and increased speed or improved accuracy in reservoir
  simulation.
\newblock {\em Multiscale Modeling \& Simulation}, 2(3):421--439, 2004.

\bibitem{aarnes2008mixed}
J{\o}rg~E Aarnes and Yalchin Efendiev.
\newblock Mixed multiscale finite element methods for stochastic porous media
  flows.
\newblock {\em SIAM Journal on Scientific Computing}, 30(5):2319--2339, 2008.

\bibitem{Arbogast_two_scale_04}
T.~Arbogast.
\newblock Analysis of a two-scale, locally conservative subgrid upscaling for
  elliptic problems.
\newblock {\em SIAM Journal on Numerical Analysis}, 42(2):576--598
  (electronic), 2004.

\bibitem{Arbogast_Boyd_06}
T.~Arbogast and K.J. Boyd.
\newblock Subgrid upscaling and mixed multiscale finite elements.
\newblock {\em SIAM J. Numer. Anal.}, 44(3):1150--1171 (electronic), 2006.

\bibitem{Arbogast_PWY_07}
T.~Arbogast, G.~Pencheva, M.F. Wheeler, and I.~Yotov.
\newblock A multiscale mortar mixed finite element method.
\newblock {\em Multiscale Model. Simul.}, 6(1):319--346, 2007.

\bibitem{arbogast2007multiscale}
Todd Arbogast, Gergina Pencheva, Mary~F Wheeler, and Ivan Yotov.
\newblock A multiscale mortar mixed finite element method.
\newblock {\em Multiscale Modeling \& Simulation}, 6(1):319--346, 2007.

\bibitem{babuvska1983generalized}
Ivo Babu{\v{s}}ka and John~E Osborn.
\newblock Generalized finite element methods: their performance and their
  relation to mixed methods.
\newblock {\em SIAM Journal on Numerical Analysis}, 20(3):510--536, 1983.

\bibitem{bush2013application}
Lawrence Bush and Victor Ginting.
\newblock On the application of the continuous galerkin finite element method
  for conservation problems.
\newblock {\em SIAM Journal on Scientific Computing}, 35(6):A2953--A2975, 2013.

\bibitem{online_mixed}
H.~Chan, E.~T. Chung, and Y.~Efendiev.
\newblock Adaptive mixed {GM}s{FEM} for flows in heterogeneous media.
\newblock {\em Numerical Mathematics: Theory, Methods and Applications},
  9(4):497--527, 2016.

\bibitem{chan2016adaptive}
Ho~Yuen Chan, Eric Chung, and Yalchin Efendiev.
\newblock Adaptive mixed gmsfem for flows in heterogeneous media.
\newblock {\em Numerical Mathematics: Theory, Methods and Applications},
  9(4):497--527, 2016.

\bibitem{chen2003mixed}
Zhiming Chen and Thomas Hou.
\newblock A mixed multiscale finite element method for elliptic problems with
  oscillating coefficients.
\newblock {\em Mathematics of Computation}, 72(242):541--576, 2003.

\bibitem{online_cg}
E.~T. Chung, Y.~Efendiev, and T.~Leung.
\newblock Residual-driven online generalized multiscale finite element methods.
\newblock {\em Journal of Computational Physics}, 302:176--190, 2015.

\bibitem{CHUNG201669}
Eric Chung, Yalchin Efendiev, and Thomas~Y. Hou.
\newblock Adaptive multiscale model reduction with generalized multiscale
  finite element methods.
\newblock {\em Journal of Computational Physics}, 320:69 -- 95, 2016.

\bibitem{chung2015mixed}
Eric~T Chung, Yalchin Efendiev, and Chak~Shing Lee.
\newblock Mixed generalized multiscale finite element methods and applications.
\newblock {\em Multiscale Modeling \& Simulation}, 13(1):338--366, 2015.

\bibitem{chung2015residual}
Eric~T Chung, Yalchin Efendiev, and Wing~Tat Leung.
\newblock Residual-driven online generalized multiscale finite element methods.
\newblock {\em Journal of Computational Physics}, 302:176--190, 2015.

\bibitem{chung2014adaptive}
Eric~T Chung, Yalchin Efendiev, and Guanglian Li.
\newblock An adaptive gmsfem for high-contrast flow problems.
\newblock {\em Journal of Computational Physics}, 273:54--76, 2014.

\bibitem{cockburn2002local}
Bernardo Cockburn, Guido Kanschat, Dominik Sch{\"o}tzau, and Christoph Schwab.
\newblock Local discontinuous galerkin methods for the stokes system.
\newblock {\em SIAM Journal on Numerical Analysis}, 40(1):319--343, 2002.

\bibitem{cortinovis2014iterative}
Davide Cortinovis and Patrick Jenny.
\newblock Iterative galerkin-enriched multiscale finite-volume method.
\newblock {\em Journal of Computational Physics}, 277:248--267, 2014.

\bibitem{du2018adaptive}
Jie Du and Eric Chung.
\newblock An adaptive staggered discontinuous galerkin method for the steady
  state convection--diffusion equation.
\newblock {\em Journal of Scientific Computing}, 77(3):1490--1518, 2018.

\bibitem{durlofsky1991numerical}
Louis~J Durlofsky.
\newblock Numerical calculation of equivalent grid block permeability tensors
  for heterogeneous porous media.
\newblock {\em Water resources research}, 27(5):699--708, 1991.

\bibitem{efendiev2002numerical}
Y~Efendiev and LJ~Durlofsky.
\newblock Numerical modeling of subgrid heterogeneity in two phase flow
  simulations.
\newblock {\em Water Resources Research}, 38(8):3--1, 2002.

\bibitem{egw10}
Y.~Efendiev, J.~Galvis, and X.H. Wu.
\newblock Multiscale finite element methods for high-contrast problems using
  local spectral basis functions.
\newblock {\em Journal of Computational Physics}, 230:937--955, 2011.

\bibitem{efendiev2009multiscale}
Y.~Efendiev and T.~Y. Hou.
\newblock {\em Multiscale finite element methods: theory and applications},
  volume~4.
\newblock Springer Science \& Business Media, 2009.

\bibitem{efendiev2000modeling}
Yalchin Efendiev, LJ~Durlofsky, and SH~Lee.
\newblock Modeling of subgrid effects in coarse-scale simulations of transport
  in heterogeneous porous media.
\newblock {\em Water Resources Research}, 36(8):2031--2041, 2000.

\bibitem{Efendiev2011ADD}
Yalchin Efendiev and Juan Galvis.
\newblock A domain decomposition preconditioner for multiscale high-contrast
  problems.
\newblock 2011.

\bibitem{efendiev2013generalized}
Yalchin Efendiev, Juan Galvis, and Thomas~Y Hou.
\newblock Generalized multiscale finite element methods ({GM}s{FEM}).
\newblock {\em Journal of Computational Physics}, 251:116--135, 2013.

\bibitem{DDDAS_upscale_2004}
V.~Ginting, R.~Ewing, Y.~Efendiev, and R.~Lazarov.
\newblock Upscaled modeling in multiphase flow applications.
\newblock {\em Comput. Appl. Math.}, 23(2-3):213--233, 2004.

\bibitem{guiraldello2020velocity}
Rafael~T Guiraldello, Roberto~F Ausas, Fabricio~S Sousa, Felipe Pereira, and
  Gustavo~C Buscaglia.
\newblock Velocity postprocessing schemes for multiscale mixed methods applied
  to contaminant transport in subsurface flows.
\newblock {\em Computational Geosciences}, pages 1--21, 2020.

\bibitem{hou1997multiscale}
Thomas~Y Hou and Xiao-Hui Wu.
\newblock A multiscale finite element method for elliptic problems in composite
  materials and porous media.
\newblock {\em Journal of computational physics}, 134(1):169--189, 1997.

\bibitem{jennylt03}
P.~Jenny, S.H. Lee, and H.~Tchelepi.
\newblock Multi-scale finite volume method for elliptic problems in subsurface
  flow simulation.
\newblock {\em Journal of Computational Physics}, 187:47--67, 2003.

\bibitem{jenny2003multi}
Patrick Jenny, SH~Lee, and Hamdi~A Tchelepi.
\newblock Multi-scale finite-volume method for elliptic problems in subsurface
  flow simulation.
\newblock {\em Journal of Computational Physics}, 187(1):47--67, 2003.

\bibitem{jiang2012some}
Lijian Jiang, J{\o}rg~E Aarnes, and Yalchin Efendiev.
\newblock Some multiscale results using limited global information for
  two-phase flow simulations.
\newblock {\em International Journal of Numerical Analysis \& Modeling}, 9(1),
  2012.

\bibitem{kim2013staggered}
Hyea~Hyun Kim, Eric~T Chung, and Chak~Shing Lee.
\newblock A staggered discontinuous galerkin method for the stokes system.
\newblock {\em SIAM Journal on Numerical Analysis}, 51(6):3327--3350, 2013.

\bibitem{kippe2008comparison}
Vegard Kippe, J{\o}rg~E Aarnes, and Knut-Andreas Lie.
\newblock A comparison of multiscale methods for elliptic problems in porous
  media flow.
\newblock {\em Computational Geosciences}, 12(3):377--398, 2008.

\bibitem{lunati2004multi}
Ivan Lunati and Patrick Jenny.
\newblock Multi-scale finite-volume method for highly heterogeneous porous
  media with shale layers.
\newblock In {\em ECMOR IX-9th European Conference on the Mathematics of Oil
  Recovery}, pages cp--9. European Association of Geoscientists \& Engineers,
  2004.

\bibitem{odsaeter2017postprocessing}
Lars~H Ods{\ae}ter, Mary~F Wheeler, Trond Kvamsdal, and Mats~G Larson.
\newblock Postprocessing of non-conservative flux for compatibility with
  transport in heterogeneous media.
\newblock {\em Computer Methods in Applied Mechanics and Engineering},
  315:799--830, 2017.

\bibitem{peszynska2005mortar}
Malgorzata Peszynska.
\newblock Mortar adaptivity in mixed methods for flow in porous media.
\newblock {\em Int. J. Numer. Anal. Model}, 2(3):241--282, 2005.

\bibitem{peszynska2002mortar}
Ma{\l}gorzata Peszy{\'n}ska, Mary~F Wheeler, and Ivan Yotov.
\newblock Mortar upscaling for multiphase flow in porous media.
\newblock {\em Computational Geosciences}, 6(1):73--100, 2002.

\bibitem{thomas2013numerical}
James~William Thomas.
\newblock {\em Numerical partial differential equations: conservation laws and
  elliptic equations}, volume~33.
\newblock Springer Science \& Business Media, 2013.

\bibitem{Wheeler_mortar_MS_12}
M.F. Wheeler, G.~Xue, and I.~Yotov.
\newblock A multiscale mortar multipoint flux mixed finite element method.
\newblock {\em ESAIM Math. Model. Numer. Anal.}, 46(4):759--796, 2012.

\bibitem{wu2002analysis}
X.~Wu, Y.~Efendiev, and T.~Y. Hou.
\newblock Analysis of upscaling absolute permeability.
\newblock {\em Discrete and Continuous Dynamical Systems Series B},
  2(2):185--204, 2002.

\bibitem{mortar_elliptic}
Y.~Yang, E.~T. Chung, and S.~Fu.
\newblock An enriched multiscale mortar space for high contrast flow problems.
\newblock {\em Commun. Comput. Phys.}, 23(2):476--499, 2018.

\bibitem{mortar_online}
Y.~Yang, E.~T. Chung, and S.~Fu.
\newblock Residual driven online mortar mixed finite element methods and
  applications.
\newblock {\em Journal of Computational and Applied Mathematics}, 34:318--333,
  2018.

\bibitem{yang2018online}
Yanfang Yang, Shubin Fu, and Eric~T Chung.
\newblock Online mixed multiscale finite element method with oversampling and
  its applications.
\newblock {\em Journal of Scientific Computing}, 82(2):31, 2020.

\bibitem{yang2019multiscale}
Yanfang Yang, Ke~Shi, and Shubin Fu.
\newblock Multiscale hybridizable discontinuous galerkin method for flow
  simulations in highly heterogeneous media.
\newblock {\em Journal of Scientific Computing}, 81(3):1712--1731, 2019.

\end{thebibliography}

\end{document}